\documentclass{amsart}

\newtheorem{theorem}{Theorem}[section]
\newtheorem*{maintheorem}{Theorem}
\newtheorem{lemma}[theorem]{Lemma}
\newtheorem{proposition}[theorem]{Proposition}
\newtheorem{corollary}[theorem]{Corollary}

\theoremstyle{definition}

\newtheorem*{acknowledgement}{Acknowledgement}

\theoremstyle{remark}

\renewcommand{\labelenumi}{(\roman{enumi})}

\DeclareFontFamily{U}{wncy}{}
\DeclareFontShape{U}{wncy}{m}{n}{<->wncyr10}{}
\DeclareSymbolFont{mcy}{U}{wncy}{m}{n}
\DeclareMathSymbol{\Sh}{\mathord}{mcy}{"58}
\newcommand\mylabel[1]{\label{#1}}

\newcommand{\NN}{\mathbb{N}}
\newcommand{\ZZ}{\mathbb{Z}}
\newcommand{\QQ}{\mathbb{Q}}
\newcommand{\RR}{\mathbb{R}}
\newcommand{\CC}{\mathbb{C}}

\newcommand{\FF}{\mathbb{F}}
\newcommand{\HH}{\mathbb{H}}
\newcommand{\PP}{\mathbb{P}}
\renewcommand{\AA}{\mathbb{A}}

\newcommand  {\shC}     {\mathcal{C}}

\newcommand  {\shN}     {\mathcal{N}}
\newcommand  {\shL}     {\mathcal{L}}

\newcommand  {\ab}      {{\text{\rm ab}}}

\newcommand  {\Aut}     {\operatorname{Aut}}

\newcommand  {\Card}    {\operatorname{Card}}

\newcommand  {\coker}   {\operatorname{coker}}
\renewcommand{\cong}    {\equiv}

\newcommand  {\End}     {\operatorname{End}}

\newcommand  {\Fr}      {\operatorname{Fr}}

\newcommand  {\Gal}     {\operatorname{Gal}}
\newcommand  {\GL}      {\operatorname{GL}}

\newcommand  {\gp}      {{\rm gp}}

\newcommand  {\Hom}     {\operatorname{Hom}}
\newcommand  {\Hilb}    {\operatorname{Hilb}}

\newcommand  {\id}      {\operatorname{id}}
\newcommand  {\im}      {\operatorname{im}}
\newcommand  {\ind}     {\operatorname{ind}}

\newcommand  {\Irr}     {{\operatorname{Irr}}}

\newcommand  {\Isom}    {\operatorname{Isom}}

\renewcommand  {\ker }  {\operatorname{ker}}

\newcommand  {\loc}     {{\rm loc}}

\newcommand  {\invlim}  {\varprojlim}

\newcommand  {\lra}     {\longrightarrow}

\newcommand  {\maxid}   {\mathfrak{m}}

\newcommand  {\NS}      {\operatorname{NS}}

\renewcommand{\O}       {\mathcal{O}}

\newcommand  {\ord}     {\operatorname{ord}}

\newcommand  {\Pic}     {\operatorname{Pic}}
\newcommand  {\PGL}     {\operatorname{PGL}}

\newcommand  {\pr}      {\operatorname{pr}}
\newcommand  {\Proj}    {\operatorname{Proj}}

\newcommand  {\quadand} {\quad\text{and}\quad}

\newcommand  {\ra}      {\rightarrow}

\newcommand  {\rank}    {\operatorname{rank}}
\newcommand  {\red}     {{\operatorname{red}}}
\newcommand  {\Reg}     {\operatorname{Reg}}
\newcommand  {\res}     {\operatorname{Res}}

\newcommand  {\Spec}    {\operatorname{Spec}}

\newcommand  {\Sym}     {\operatorname{Sym}}

\newcommand  {\Tr}      {\operatorname{Tr}}
\newcommand  {\Tor}     {\operatorname{Tor}}

\newcommand  {\val}     {\operatorname{val}}

\def\mydate{\number\day\space\ifcase\month \or January\or February\or March\or 
April\or May\or June\or July\or
August\or September\or October\or November\or December\fi \space\number\year}

\DeclareFontFamily{U}{wncy}{}
\DeclareFontShape{U}{wncy}{m}{n}{<->wncyr10}{}
\DeclareSymbolFont{mcy}{U}{wncy}{m}{n}
\DeclareMathSymbol{\Sh}{\mathord}{mcy}{"58}
\usepackage{amscd,amssymb,amsmath, xypic, graphics}

\begin{document}

\title[Wild quotient singularities]
      {Wild quotient surface singularities  whose dual graphs are not star-shaped}

\author[Hiroyuki Ito]{Hiroyuki Ito}
\address{Department of Mathematics,
Faculty of Science and Technology,
Tokyo University of Science,
2641 Yamazaki, Noda, Chiba, 278-8510, Japan}
\curraddr{}
\email{ito\_hiroyuki@ma.noda.tus.ac.jp}

\author[Stefan Schr\"oer]{Stefan Schr\"oer}
\address{Mathematisches Institut, Heinrich-Heine-Universit\"at,
40204 D\"usseldorf, Germany}
\curraddr{}
\email{schroeer@math.uni-duesseldorf.de}

%Singularities, Singularities of surfaces

\dedicatory{Revised version, 26 August 2014}

\begin{abstract}
We obtain  results that  answer  certain questions of Lorenzini
on wild quotient singularities in dimension two:
Using Kato's theory of log structures and log regularity, we prove that the dual graph of exceptional curves 
on the resolution of singularities contains at least one node. 
Furthermore, we show that diagonal quotients for Hermitian curves by analogues of Heisenberg groups
lead to examples of wild quotient singularities  where the dual graph contains at least two nodes.
\end{abstract}

\maketitle
\tableofcontents
\renewcommand{\labelenumi}{(\roman{enumi})}

%===========================================================
\section*{Introduction}

Quotient singularities arise in algebraic geometry and commutative algebra
whenever   finite groups act  on  complex spaces,   schemes, or rings.
Roughly speaking, a \emph{quotient singularity} is a local   ring $R=A^G$
arising as the ring of invariants    inside a regular local ring with respect to a finite
group of automorphisms. It is called \emph{wild} if the characteristic of the
residue field $k=R/\maxid_R$ is positive and divides the order of the group $G$, and \emph{tame} otherwise.
Although  quotient singularities play a crucial role in various
areas (toric varieties, McKay correspondence, stable reduction,
singularity theory,  surfaces of general type etc.), 
only few concrete examples and general results seem  to be known for the wild case.

Lorenzini initiated a systematic study of \emph{wild quotient singularities in dimension two} in
a series of papers (\cite{Lorenzini 2013}, \cite{Lorenzini 2014}, and  \cite{Lorenzini 2011c}).
The goal of this article is to further investigate such singularities.
Let $R$ be a   local noetherian ring that is normal and 2-dimensional. For simplicity we 
suppose that the ring is complete, with separably closed residue field.
Attached to this situation is the so-called \emph{dual graph}. This is the graph
whose vertices correspond to the irreducible components   of the
exceptional divisor $E\subset X$ on the minimal good resolution of singularities $X\ra\Spec(R)$,
with edges indicating intersections.
It is well-known that the dual graph of a quotient singularity is a \emph{tree}, and that
the irreducible components are copies of $\PP^1$.

Brieskorn's characteristic zero classification of quotient singularities \cite{Brieskorn 1968} 
reveals that the number of nodes in the  dual graph is either zero or one.
In stark contrast, little  seems to be known about the structure of dual graphs for wild quotient singularities.
In this context the term \emph{node} refers to vertices  of valency $\geq 3$,
that is, having at least three neighbors.

Lorenzini   recently posed a list of striking questions on such singularities \cite{Lorenzini 2011d}.
Among other things, he  asks whether there are wild quotient singularities whose \emph{dual graphs have no node},
or some whose \emph{dual graphs  have at least two nodes},
and whether the dual graphs occurring \emph{over fields are the same as those occurring
in mixed characteristics}.  Our first main result answers the first question to the negative:

\begin{maintheorem} [see Corollary \ref{one node}]
Dual graphs of   wild quotient singularities in dimension two contain at least one node.
\end{maintheorem}

We prove this by computing the \emph{local fundamental group} $\pi_1^\loc(R)$ of    2-dimensional singularities 
whose exceptional divisor is a chain of projective lines. Such   are called \emph{Hirzebruch--Jung singularities}.
Note that we do not suppose that $R$ contains a field.
It turns out that   the local fundamental group is cyclic of order prime to the characteristic
exponent of the residue field. This generalizes Artin's observation on rational double points over fields
\cite{Artin 1977}. Our approach hinges on Kato's theory of \emph{log structures} \cite{Kato 1989},
and in particular his results on \emph{log regularity} \cite{Kato 1994}, a notion that generalizes toroidal embeddings
to   mixed characteristics. We refer to Ogus's lecture notes \cite{Ogus 2012} for an in-depth exposition  of this theory.

The second part of  this paper deals with the construction of 
wild quotient singularities having at least two nodes. Here we work over a ground field $k$ of characteristic
$p>0$. The idea is rather simple: Take a curve $C$ with unusually large automorphism group,
and look at the quotient $G\backslash(C\times C)$ of the selfproduct by the diagonal action of the
Sylow-$p$-subgroup $G\subset\Aut(C)$.
Quotients of the form $G\backslash(C\times C')$ were already studied by Lorenzini \cite{Lorenzini 2011c},
who treated the case that one factor is ordinary, and ourselves \cite{Ito; Schroeer 2012} in the case
that both factors are certain Artin--Schreier curves, chosen in such a way that local computations
become feasible.

Here we study the  case where 
$$
C:\quad y^q-y=x^{q+1}
$$ 
are so-called \emph{Hermitian curves}, for prime powers $q=p^m$.
By results of Stichtenoth, the number of automorphisms way exceeds   the Hurwitz bound.
In fact, these curves have the largest possible number of automorphisms in relation
to the genus, even by characteristic $p$ standards (\cite{Stichtenoth 1973a} and \cite{Stichtenoth 1973b}). 
The Sylow-$p$-subgroup $G\subset\Aut(C)$ has order $\ord(G)=q^3$.
It turns out that $G$ is a so-called \emph{special $p$-group}.
These are certain nonabelian groups analogous to Heisenberg groups, which play
a   role in the classification of the finite simple groups. All in all,
$G\backslash(C\times C)$ is a very natural candidate to look for non-star-shaped dual graphs.
Our second main result tells us that this indeed something   happens:

\begin{maintheorem} [see Theorem \ref{two nodes}]
The dual graph for the minimal resolution of singularities for the wild quotient
singularity on $G\backslash(C\times C)$ contains at least two nodes.
\end{maintheorem}

Note that our arguments are   indirect, such that we gain no control over the precise form
of the dual graph nor its exact number of nodes, let alone the formal equations for the singularity.
In some sense, this answers Lorenzini's second question in a positive way. 
He had, however, only the case $G=\ZZ/p\ZZ$ in mind. Note  also that he had already constructed
an example of a wild  quotient singularity in mixed characteristics over
the 2-adic numbers with two nodes \cite{Lorenzini 2013}.

Finally, we analyze the induced group action on the \'etale cohomology $H^1(C,\QQ_l)$.
Note that there are many formulas expressing   cohomology   as representations
for tame actions, for example due to  
Chevalley and Weil \cite{Chevalley; Weil 1934}, Ellingsrud and L\o nsted \cite{Ellingsrud; Lonsted 1980},  K\"ock \cite{Kock 2005},
and others, but all breaks down in the wild case. Rather,  we have to apply brute force,
and completely determine the irreducible representations of our special $p$-groups
over various fields, and then single out the given representation via the Lefschetz Trace Formula.
Our third main result is an explicit description of the  cohomology as representation.
It might be stated as follows:

\begin{maintheorem}[see Theorem \ref{cohomology representation}]
Let $l\neq p$ be a prime different from the characteristic $p$.
In case $p\geq 3$ we assume that $p$ does not divide $l-1$.
Then the $G$-representation $H^1(C,\QQ_l)$ is the direct sum over a basic set of irreducible
$G$-representations over $\QQ_l$ that do not factor over the abelianization $G^\ab$.
\end{maintheorem}

See Sections \ref{Representations} for an explicit description of irreducible representations over various fields.
It turns out that under the assumptions of the preceding Theorem, the cohomology $H^1(C,\QQ_l)$
is irreducible as a representation of the full automorphism group $\Aut(C)$, which appears to be rather peculiar.
Furthermore, knowing cohomology as representation   allows us to compute the
Chern numbers for the relative minimal model of $G\backslash(C\times C)$.

The paper is organized as follows. Section \ref{Wild quotient singularities} contains
general facts on wild quotients. In Section \ref{Dual graphs without nodes} we give
an abstract criterion to express certain local fundamental group as Galois groups
and state our result on the local fundamental group of Hirzebruch--Jung singularities.
We prove this result in Section \ref{Local fundamental groups}, using Kato's theory
of log schemes and log regularity. A technical step, which is elementary but somewhat lengthy, 
is completed in Section \ref{Cotangent spaces}.
In Section \ref{Hermitian curves} we turn  to Hermitian curves and their automorphism groups.
In Section \ref{Two nodes} we state that the corresponding diagonal quotients yield wild quotient
singularities with at least two nodes. The proof occupies Sections \ref{Generic fiber} and \ref{Singular fiber}.
We return to Hermitian curves in Section \ref{Higher Ramification}, where the ramification groups
and Swan conductors are computed. In Section \ref{Representations} we examine the extraspecial $p$-groups acting on Hermitian curves
and make a comprehensive investigation of their   representation theory. This is applied in Section \ref{Cohomology as},
in which the first $l$-adic cohomology of Hermitian curves is computed as a representation.
In Section \ref{Curves with}, we observe
that such curves are supersingular. Section \ref{Diagonal quotients} contains a general investigation of diagonal quotients
and their Picard numbers in terms of cohomology as representation.
We apply this in the final Section \ref{Chern invariants} to compute the Chern numbers of the minimal smooth
model of the diagonal quotient for   Hermitian curves.

\begin{acknowledgement}
The first author would like to thank the Mathematisches Institut of the Heinrich-Heine-Universit\"{a}t D\"{u}sseldorf, 
where this work has progressed, for its warm hospitality. Research of the first author was partially supported by 
Grand-in-Aid for Scientific Research (C) 24540051, MEXT. 

We also would like to thank the referee for careful reading and suggestions to improve the paper.
\end{acknowledgement}

%===========================================================
\section{Generalities on wild quotient singularities}
\mylabel{Wild quotient singularities} 

Let   $G$ be a finite group acting on a regular local noetherian ring $A$.
Then the invariant ring $R=A^G$ is a normal local noetherian ring, and the ring extension $R\subset A$ is
finite. If the induced $G$-action on $X=\Spec(A)$
is free in codimension one, the ring $R$ and its spectrum $Y=\Spec(R)$
are called \emph{quotient singularity}.  Note that we do not assume that our rings contain   fields.

At various places it will be convenient to use the  \emph{characteristic exponent}   of a  field; recall that this
is 1 if the field has characteristic zero, and   equal to the characteristic of the field otherwise.
A quotient singularity $R$ is called \emph{tame} if the order of $G$ is prime to the
characteristic  exponent $p\geq 1$ of the residue field $k=A/\maxid_A$. The  situation becomes particularly simple if 
$A=\O_{X,x}^\wedge$ is a complete local ring for a point $x\in X$ on a complex manifold:
By results going back to H.\ Cartan \cite{Cartan 1957},  the group $G$ acts linearly on a suitable regular system of parameters. 
If $G$ is abelian, one may diagonalize simultaneously, and the situation is easily described in terms of toric varieties.

A quotient singularity $R$ is called \emph{wild} if the order of $G$ is not prime to the characteristic exponent 
of the residue field. Many properties of tame quotient singularities do not hold true: 
First of all, it is unclear   whether or not \emph{resolutions of singularities} exists in higher dimensions.
Next, wild quotient singularities are usually not \emph{Cohen--Macaulay}.
This was first observed by Fossum and Griffith \cite{Fossum; Griffith 1975}, and studied in more detail by
Ellingsrud and Skjelbred \cite{Ellingsrud; Skjelbred 1980}, for groups $G$ acting on polynomial rings
by permutation of the indeterminates. 

Another problem is that  the \emph{Grauert--Riemenschneider Vanishing Theorem} does not necessarily hold. 
To give an explicit  example, let $S$ be a smooth surface over a ground field $k$, 
and consider the symmetric product $\Sym^n(S)=(S\times\ldots\times S)/G$, where $G$ is the
symmetric group on $n$ letters acting by permutation of the factors.
The Hilbert--Chow morphisms gives a resolution of singularities 
$$
f:X=\Hilb^n(S)\lra\Sym^n(S)=Y,
$$
which is a crepant resolution. If the characteristic $p$ is positive with $p\leq n$, then 
$R^if_*(\omega_{X/Y})=R^if_*(\O_X)\neq 0$, because otherwise the augments
in   \cite{Kempf et al. 1973}, page 49--51 show  that $Y$ would be Cohen--Macaulay. But this  is not the case
by the work of Ellingsrud and Skjelbred.  

In this paper we are mainly interested in \emph{wild quotient singularities in dimension two}.
Let $R$ be an  arbitrary local noetherian ring that is normal. To simplify,
we also assume that the residue field is separably closed and $R$ is complete.  
The latter ensures that there is a resolution of singularities
$f:X\ra\Spec(R)$ without further assumptions, according to Lipman  \cite{Lipman 1978}.
Let $E\subset X$ be the exceptional divisor, viewed as a reduced closed subscheme, 
and $E_1,\ldots, E_r$ be its irreducible components.
We shall consider the \emph{intersection matrix} $(E_i\cdot E_j)_{1\leq i,j\leq r}$ 
and the \emph{dual graph}, which is the graph whose vertices corresponds to the irreducible components
$E_i$. Two   vertices are joined by an edge if $E_i\cap E_j\neq\emptyset$.
Here we tacitly assume that $E$ has simple normal crossings.
It is well-known that   dual graphs for quotient singularities $R=A^G$ in dimension two are necessarily   \emph{trees},
and that all irreducible components are $E_i\simeq \PP^1$ 
(see \cite{Cutkosky; Srinivasan 1993}, Theorem 1 or \cite{Lorenzini 2013}, Theorem 2.8).

Brieskorn \cite{Brieskorn 1968} gave a complete classification of quotient singularities in characteristic zero,
which builds on the ADE-classification of the rational double points. 
From this one sees that dual graphs of tame quotient singularities are always \emph{star-shaped},
that is, contain  at most one node. Absence of nodes means that the group $G$ is abelian.
Following Lorenzini, we call  \emph{nodes} the vertices of valency $\geq 3$, that is,  with at least three neighboring vertices.

In the wild case, the   dual graphs may become rather complicated.
Shioda \cite{Shioda 1974} and later Katsura \cite{Katsura 1978} analyzed the Kummer surface $\left\{\pm 1\right\}\backslash A$
in characteristic 2, and observed   certain minimally elliptic singularities.
In \cite{Ito; Schroeer 2012}, we studied quotients $G\backslash(C\times C')$, where the curves $C,C'$ are Artin--Schreier curves
of the form $y^q-y=f(x)$, where the right hand side is a polynomial of degree $q-1$. It turned out that
the corresponding wild quotient singularity is star-shaped, with $q+1$ chains of length $q-1$ sprouting from
the node. Moreover, the geometric genus of this wild quotient singularity grows at least quadratically in $q$. Here $q$ denotes
an arbitrary $p$-power.
Among other things, Lorenzini \cite{Lorenzini 2011d}  asks whether there are wild quotient singularities whose \emph{dual graph has no node},
or some whose \emph{dual graph  has at least two nodes},
and whether the dual graphs occurring \emph{over fields are the same as those occurring
in mixed characteristics}.

%===========================================================
\section{Dual graphs without nodes}
\mylabel{Dual graphs without nodes} 

Let $R$ be a complete local noetherian ring 
that is normal and 2-dimen\-sio\-nal, and whose residue field $k=R/\maxid_R$ is separably closed.
Note that we do not assume that $R$ contains a field. 
Let $f:X\ra\Spec(R)$ be the minimal  resolution of singularities, with exceptional divisor $E\subset X$.
Let $E_1,\ldots,E_r\subset E$ be the integral components, and 
$(E_i\cdot E_j)_{1\leq i,j\leq r}$ be the intersection matrix, which is negative-definite.
We say that $R$ is a \emph{Hirzebruch--Jung singularity} if all $E_i\simeq\PP^1_k$, and 
for some suitable ordering of the irreducible components,
the intersection matrix has the form
\begin{equation}
\label{intersection matrix}
(E_i\cdot E_j)_{1\leq i,j\leq r} =
\begin{pmatrix}
-s_1 & 1\\
1    & -s_2 & 1 \\
     &  1   & -s_3   & \ddots\\
     &      & \ddots & \ddots & 1\\
     &      &        & 1      & -s_r
\end{pmatrix}
\end{equation}
for some integers $s_i\geq 2$. Whence the dual graph attached to the singularity is
a chain:

%-------------------------------------------------------------
\vspace{2em}
\centerline{\includegraphics{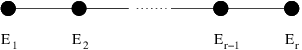} }
\vspace{1em}
%\centerline{Figure \stepcounter{figure}\arabic{figure}: Dual graph of the exceptional divisor $E\subset X$.}
\vspace{1em}
%-------------------------------------------------------------
We may assign a reduced fraction $m/b\in \QQ_{>1}$ to the Hirzebruch--Jung singularity, defined by the continued fraction
\begin{equation}
\label{continued fraction}
\frac{m}{b} = 
s_1- \cfrac {1}{s_2-
    \cfrac{1}{\ddots -
     \cfrac{1}{s_r}}}.        
\end{equation}
Writing the exceptional curves in   reverse order, one obtains the fraction $m/b'$ where $0<b'<m$ is the
unique integer with $bb'\cong 1$ modulo $m$.
We say that $R$ is a \emph{Hirzebruch--Jung singularity of type} $m/b$. Note that 
$R$ is also a Hirzebruch--Jung singularity of type $m/b'$.
In the special case $b=m-1$ we have   $s_1=\ldots=s_r=2$ and $r=m-1$. In other words, the
Hirzebruch--Jung singularities of type $m/(m-1)$ are precisely the \emph{rational double points}
of type $A_{m-1}$. In the special case $r=1$ we have $m=s_1$ and $b=1$,
and the Hirzebruch--Jung singularities of type $m/1$ are obtained by contracting
a curve $E=E_i\simeq\PP^1$ of self-intersection $-m$.

Recall that the \emph{local fundamental group} $\pi_1^\loc(R)$ of a complete local noetherian domain $R$
is the fundamental group of the open subscheme $\Reg(R)\subset\Spec(R)$, say with respect
to a chosen separable closure of the field of fractions.
If $R$ is a normal and 2-dimensional, it is possible to compute the \emph{maximal  quotient} of 
$\pi_1^\loc(R)$ that is prime to the characteristic exponent of the residue field,
in terms of the intersection matrix of the exceptional divisor, as explained
by Cutkosky and Srinivasan \cite{Cutkosky; Srinivasan 1993}, Theorem 3.
Our first main result is a computation of the \emph{full} local fundamental group:

\begin{theorem}
\mylabel{fundamental group}
Let $R$ be a Hirzebruch--Jung singularity of type $m/b$, and $p\geq 1$ the characteristic exponent
of the residue field $k=R/\maxid_R$.
Then the local fundamental group $\pi_1^\loc(R)$ is isomorphic to the   prime-to-$p$
part of the cyclic group $\ZZ/m\ZZ$.
\end{theorem}

Of course, this is well-known in the case that $R$ contains a field of characteristic zero.
For rational double points of type $A_{m-1}$, with $R$ containing a field of characteristic $p>0$, 
it was remarked by Artin
(see \cite{Artin 1977}, Section 2). The really new part in our result 
seems to be the case of mixed characteristics. The proof will stretch over the next two sections.

The theorem implies that the local fundamental group of a Hirzebruch--Jung singularity   contains no   elements of order $p$.
Whence the following immediate consequence for wild quotient singularities,
which answers a question of Lorenzini  to the negative (\cite{Lorenzini 2011d}, Question 1.1 (a)):

\begin{corollary}
\mylabel{one node}
The dual graph attached to  a wild quotient singularity in dimension two  contains at least one node.
\end{corollary}

To compute the local fundamental group of a Hirzebruch--Jung singularity,
we shall use the following   fact,
which   already appeared in a special form in \cite{Artin 1977}, Corollary 1.2. We state it in a
rather general way: 
Let $A\subset A'\subset A''$ be three complete local noetherian normal domains, with respective fields
of fractions $K\subset K'\subset K''$. Assume that the extension $A\subset A'$ and $A'\subset A''$ are
finite.

\begin{proposition}
\mylabel{galois group}
In the preceding situation, suppose the following:
\begin{enumerate}
\item The local ring $A$ has separably closed residue field.
\item The local ring $A''$ is regular.
\item The finite field extension $K\subset K'$ is Galois.
\item The ring extension $A\subset A'$ is \'etale in codimension one.
\item The ring extension $A'\subset A''$ is totally ramified at some prime $\mathfrak{q}\subset A'$ of height one.
\end{enumerate}
Then there is a canonical identification $\pi_1^\loc(A)=\Gal(K'/K)$, where the former is taken
with respect to some separable closure of $K'$.
\end{proposition}

% Suppose we have a commutative diagram
% $$
% \begin{CD}
% A @>>> A' @>>> A''\\
% @VVV @VVV\\
% K @>>> K'
% \end{CD}
% $$
% where $A,A',A''$ are complete local noetherian rings that are normal,
% $A\subset A'\subset A''$ are finite ring extensions,
% and $A\subset K$ and $A'\subset K'$ are the fields of fractions.
% 
% \begin{lemma}
% \mylabel{galois group}
% Suppose that the ring $A$ has separably closed residue field,
% that the ring $A''$ is regular, that the   field extension $K\subset K'$ is Galois,
% that $A\subset A'$ is \'etale in codimension one, and that $A'\subset A''$
% is totally ramified at some prime $\mathfrak{q}\subset A'$ of height one.
% Then there is a canonical identification $\pi_1^\loc(A)=\Gal(K'/K)$.
% \end{lemma}

\proof
This identification arises as follows:
Let $Y\subset\Spec(A)$   be the regular locus, which is an
open subscheme  by Nagata's Theorem (\cite{EGA IVb}, Theorem 6.12.7), and $Y'\subset \Spec(A')$ be its preimage.
Consider the \emph{Galois categories} $\shC$ of all finite \'etale morphisms $\tilde{Y}\ra Y$,
and $\shC'$ of all finite left $\Gal(K'/K)$-sets (compare \cite{SGA 1}, Expos\'e V). We have a  functor
$$
\Psi:\shC\lra\shC',\quad \tilde{Y}\longmapsto \tilde{Y}(K'),
$$
where $\tilde{Y}(K')$ denotes the set of $K'$-rational points of the finite $K$-scheme $\tilde{Y}_K=\tilde{Y}\times_Y\Spec(K)$.
The task now is to show that this functor is an equivalence of categories.
Indeed, from this it follows that the set-valued functor that forgets the Galois action yields compatible fiber functors
on $\shC$ and $\shC'$, such that $\pi_1^\loc(A)$  and $\Gal(K'/K)$, which  by definition are the automorphism groups of
these fiber functors, become  equal.

It is easy to see that the functor $\Psi$ is faithful:
The restriction functor $\tilde{Y}\mapsto\tilde{Y}_K$ is faithful, because $\O_Y$ is torsion free,
and the functor $\tilde{Y}_K\mapsto \tilde{Y}(K)$ is faithful by Galois descent.

We next check that $\Psi$ is essentially surjective:
Set $G=\Gal(K'/K)$.
Since $A\subset A'$ is \'etale in codimension one and both rings are normal, 
the projection $Y'\ra Y$ must be a $G$-torsor with respect to the canonical right $G$-action on $Y'$,
by the Zariski--Nagata Purity Theorem (\cite{SGA 2}, Expos\'e X, Theorem 3.4).
For every subgroup $H\subset G$, the scheme $\tilde{Y}=Y'/H$ is \'etale over $Y$,
with $\tilde{Y}(K')=G/H$. Using that every finite left $G$-set is
isomorphic to a disjoint union of sets of the form $G/H$, we infer that $\Psi$ is essentially surjective.

It remains to verify that $\Psi$ is full.
For this, it suffices to show that
every finite \'etale $g:\tilde{Y}\ra Y$ acquires a section
over $K'$, because then the scheme $\tilde{Y}_{K'}$ becomes the disjoint union of copies of $\Spec(K')$ indexed
by the finite set $\tilde{Y}(K')$. The latter ensures that any two morphisms $h,h':\tilde{Y}\ra\tilde{Z}$ between finite
\'etale $Y$-scheme that induce the same map on $K'$-rational points must be equal.

By We first extend the quasifinite morphism
$$
\tilde{Y}''=\tilde{Y}\times_Y\Spec(A'')\ra\Spec(A'')
$$
to a finite morphism $\Spec(\tilde{A}'')\ra\Spec(A'')$ with $\tilde{A}''$ normal, using Zariski's Main Theorem
\cite{EGA IVd}, Corollary 18.12.13.
Then $A''\subset \tilde{A}''$  is \'etale in codimension one by construction.
Since $A''$ is regular, the  finite ring extension must be be \'etale everywhere, again by the Zariski--Nagata Purity
Theorem. Since  $A''$ is strictly henselian, the \'etale extension $A''\subset\tilde{A}''$ admits
a retraction, by \cite{EGA IVd}, Proposition 18.8.1. Choosing such a retraction,
we obtain a factorization $Y''\ra\tilde{Y'}\ra Y'$, where $Y''\subset\Spec(A'')$
is the preimage of $Y$, and $\tilde{Y}'=\tilde{Y}\times_Y Y'$.
We conclude that $\widetilde{Y}'\ra Y'$ 
is totally ramified on some connected component of $\tilde{Y}'$ over the point $y'\in Y'$ corresponding to $\mathfrak{q}\subset A'$.
Being finite and \'etale,   $\tilde{Y}'\ra Y'$ must have a section, in particularly $\tilde{Y}(K')\neq\emptyset$.     
\qed

%===========================================================
\section{Local fundamental groups}
\mylabel{Local fundamental groups} 

Let $R$ be a Hirzebruch--Jung singularity of type $m/b$, as in the preceding section,
with residue field $k=R/\maxid_R$.
The goal of this section is to compute the local fundamental group $\pi_1^\loc(R)$, and thus to prove Theorem \ref{fundamental group}.

Let $E=E_1+\ldots+E_r$ be the reduced exceptional divisor on the minimal good
resolution of singularities $f:X\ra\Spec(R)$.
Recall that the \emph{fundamental cycle} $Z\subset X$ is the smallest
exceptional divisor containing $E$ with $(Z\cdot E_i)\leq 0$ for all $1\leq i\leq r$.

\begin{proposition}
\mylabel{rational singularity}
Our Hirzebruch--Jung singularity $R$ has fundamental cycle $Z=E$,
and is a rational singularity. Also, $E\subset X$ is the schematic closed
fiber of the resolution $f:X\ra\Spec(R)$.
\end{proposition}

\proof
Since $(E\cdot E_i)\leq 2-s_i\leq 0$
we must have $Z=E$. Clearly, $H^1(Z,\O_Z)=0$.
Whence the singularity is rational by Artin's Criterion 
\cite{Artin 1966}, Theorem 3. The final statement is contained in  loc.\ cit., Theorem 4.
\qed

\medskip
The main idea now is to use certain \emph{log structures} on $\Spec(R)$, in the sense of Kato \cite{Kato 1989}.
To this end, choose $k$-valued points $x\in E_1\smallsetminus E_2$ and $x'\in E_r\smallsetminus E_{r-1}$.
We may extend these points to effective Cartier divisors $D,D'\subset X$ with $D\cap E=\left\{x\right\}$ and
$D'\cap E=\left\{x'\right\}$, because $R$ is henselian (compare \cite{EGA IVd}, Proposition 21.9.11).
Their images
$$
C=f(D) \quadand C'=f(D')
$$
are Weil divisors on $\Spec(R)$.
Consider the   submonoid
$$
M = \left\{g\in R\mid \text{  $\Spec(R/gR)$ is supported by $C\cup C'$}\right\}
$$
inside the multiplicative monoid $R$.
It comes with a homomorphism of monoids
$$
\nu:M\lra \NN\oplus \NN,\quad g\longmapsto (\val_C(g),\val_{C'}(g)),
$$
sending $g\in M$ to the values of the valuations corresponding to the generic points
of $C,C'$. Let $P\subset\NN\oplus\NN$ be its image. Then we have a sequence of monoids
$$
1\lra R^\times \lra M \lra P \lra 0.
$$
This sequence is exact, in the sense that the monoid $P$ is the quotient of the monoid $M$
by the congruence relation $\left\{(g,g')\mid \text{$g=ug'$ for some $u\in R^\times$}\right\}$.

We may describe $P$ in terms of the intersection matrix as follows. Let
$\Phi$ be the \emph{discriminant group} of the intersection matrix $(E_i\cdot E_j)_{1\leq i,j\leq r}$,
that is, the cokernel of the map $\ZZ^{\oplus r}\ra\ZZ^{\oplus r}$ defined by the intersection matrix.
Its structure is well-known:

\begin{proposition}
The discriminant group $\Phi$ is cyclic of order $m$.
\end{proposition}

\proof
Recall that the reduced fraction $m/b$ is defined as  the continued fraction (\ref{continued fraction}).
Glancing at the intersection matrix (\ref{intersection matrix}), we see that 
it contains as submatrix an $(r-1)\times(r-1)$ identity matrix, so $\Phi$ must be cyclic. By induction on $r\geq 1$, one
easily infers 
$$
m=(-1)^r\det(E_i\cdot E_j)_{1\leq i,j\leq r}\quadand
b=(-1)^{r-1}\det(E_i\cdot E_j)_{2\leq i,j\leq r},
$$
whence $\Phi$ has order $m$.
\qed

\medskip
Consider the homomorphism of groups $\varphi:\ZZ\oplus\ZZ\ra\Phi$ sending $(n,n')$ to the class
of the $r$-tuple $(n,0,\ldots,0,n')$.

\begin{proposition}
\mylabel{description P}
We have $P=\ker(\varphi)\cap(\NN\oplus\NN)$ inside $\ZZ\oplus\ZZ$.
\end{proposition}

\proof
Let $n,n'\geq 0$ be natural numbers.
Suppose $(n,n')\in P$, such that $nC+n'C'$ is Cartier. Its preimage $nD+\sum\lambda_jE_j+n'D'$
is a Cartier divisor on $X$ that is numerically trivial on $E$.
Whence $(nD+n'D')\cdot E_i= -\sum\lambda_jE_j\cdot E_i$ for all $1\leq i\leq r$, 
 thus $\varphi(n,n')=0$.

Conversely, if the tuple of intersection numbers $(n,0,\ldots,0,n')$ is
of the form $(-\sum\lambda_jE_j\cdot E_i)_{1\leq i\leq r}$,
then the invertible sheaf $\shL=\O_X(nD+\sum\lambda_jE_j+n'D')$ is numerically trivial on $E$. It must be
trivial on the formal neighborhood of $E$, since the singularity $R$ is rational,
thus $f_*(\shL)$ remains invertible. Hence this is also the reflexive sheaf of rank one
attached to the Weil divisor $nC+n'C'$, therefore $(n,n')\in P$.
\qed

\medskip
This ensures that our monoid $P$ has all sorts of good properties.
Throughout, all monoids are assumed to be commutative.
Recall that a monoid $Q$ is called \emph{fine} if it is finitely generated and integral.
It is called \emph{saturated} if it is integral, and the inclusion $Q\subset Q^\gp$ into its groupification
has the property that whenever $a\in Q^\gp$ has a multiple in $Q$, then it is already contained in $Q$.
A fine saturated monoid is called \emph{fs}. A monoid $Q$ is \emph{sharp} if $Q^\times=\left\{1\right\}$.
We refer to Kato's original article \cite{Kato 1989} or the in-depth manuscript of Ogus \cite{Ogus 2012} for 
all standard fact on monoids.

\begin{corollary}
\mylabel{properties P}
The monoid $P$ is sharp and fs, and the surjective homomorphism of monoids $M\ra P$ has a section.
Moreover, $(\ZZ\oplus\ZZ)/P^\gp$ is  a cyclic group of order $m$.
\end{corollary}

\proof
By definition, $P=M/R^\times$ is sharp. In light of  Proposition \ref{description P}, the monoid $P$ can be
described as the set of integral solutions of finitely many inequalities with integral coefficients,
and is thus fs (for example, \cite{Grillet 2001}, Proposition 8.5).
It is easy to see that $P^\gp=\ker(\varphi)$, such that $(\ZZ\oplus\ZZ)/P^\gp=\Phi$ is cyclic of order $m$.
A section exists because the abelian group $P^\gp$ is free, thus admits only trivial extensions.
\qed

\medskip
This has a geometric interpretation in terms of toric varieties:
Let $F$ be an arbitrary ring. The inclusion $P^\gp\subset\ZZ\oplus\ZZ$
of free abelian groups of rank $2$
corresponds to an isogeny $\Spec F[\ZZ\oplus\ZZ]\ra\Spec F[P^\gp]$
between families of 2-dimensional tori over $F$. Its kernel is
the diagonalizable groups scheme $H=\Spec F[N]$,
which represents the functor on $F$-algebras $A\mapsto\Hom(N,A^\times)$,
as explained in \cite{SGA 3b}, Expose VIII.
This isogeny extends to a toric morphism $\Spec F[\NN\oplus\NN]\ra\Spec F[P]$
between families of affine toric surfaces over $F$. In fact, this toric morphism
is the quotient morphisms for the canonical $H$-action on $\AA^2_F=\Spec F[\NN\oplus\NN]$.
Let us regard $\Spec F$ as the closed subscheme of $\Spec F[P]$ corresponding to the
ideal $(P\smallsetminus P^\times)\subset F[P]$.
The following seems to be well-known:

\begin{proposition}
\mylabel{free action}
The quotient morphism $\Spec F[\NN\oplus\NN]\ra\Spec F[P]$ is a $H$-torsor
over the complement of $\Spec(F)\subset\Spec(F[P])$.
\end{proposition}

\proof
It suffices to check this on geometric fibers, so we may assume that
$F$ is an algebraically closed field.
The $H$-action commutes with the torus action, whence the torus orbits
are $H$-invariant. It suffices to check that $H$ acts freely 
on each of the the three non-closed torus orbits. This
is clear for the dense torus orbit. 
For the remaining two non-closed torus orbits, it suffices to check
that the subgroups $\ZZ\times 0$ and $0\times\ZZ$ surject onto $N$
under the projection $\ZZ^2\ra N$. This is independent of the ground field $F$,
so it suffices to treat the case $F=\CC$. Then the statement
follows from the theory of complex toric surface singularities:
According to \cite{Barth; Peters; Van de Ven 1984}, Theorem 5.1,
all complex Hirzebruch--Jung singularities of type $m/b$ are formally
isomorphic, and in particular isomorphic to the toric surface singularity
coming from the convex rational polyhedral cone $\RR_{\geq 0}(1,0)+\RR_{\geq 0}(m,-b)\subset\RR^2$,
see \cite{Cox; Little; Schenck 2011}, Theorem 10.2.3. It then follows
from loc.\ cit., Proposition 10.1.2 that the $H$-action on
$\Spec \CC[\NN\oplus\NN]$ is free outside the origin.
\qed

\medskip
To proceed, choose a section of $M\ra P$, and write $M=R^\times\oplus P$. In this way, we regard $P$ as 
a submonoid of $R$, and   revert to multiplicative notation.
Now recall that every monoid $Q$ has a unique maximal ideal,
namely the complement of the units $Q\smallsetminus Q^\times$.
For sharp monoids, this is just $Q\smallsetminus\left\{1\right\}$.
We  come to a key observation:

\begin{lemma}
\mylabel{maximal ideal}
For every choosen section of $M\ra P$, the maximal ideal  $P\smallsetminus P^\times$ of the monoid   generates the maximal
ideal $\maxid_R\subset R$ of the ring.
\end{lemma}

The proof is elementary but somewhat lengthy, so we defer it into the next section.
To proceed, we have to distinguish two cases:
If $R$ contains a field,   $W\subset R$ denotes a \emph{field of representatives},
that is, a subfield so that the projection $W\ra R/\maxid_R$ is bijective.
If $R$ does not contain a field,   choose a \emph{Cohen subring} $W\subset R$.
In our situation, this is a complete discrete valuation ring so that $p\in W$
is a uniformizer, and $W/pW\ra R/\maxid_R$ is bijective
(compare \cite{AC 8-9}, Chapter IX, \S2-3).

In any case, Lemma \ref{maximal ideal} ensures that the   map 
$$
W[[P]]\lra R
$$ 
is surjective. Note that this map is defined with the help of the choosen section
of $M\ra P$. Here $W[[P]]$ is the ring of formal power sequences $\sum_{g\in P} w_gg$ with coefficients
$w_g\in W$. The  
ring structure comes from the interpretation $W[[P]]=\invlim W[P]/I^nW[P]$ as $I$-adic
completion of the monoid ring, where $I\subset W[P]$ is the ideal generated by $P\smallsetminus P^\times$.
In Kato's terminology \cite{Kato 1994}, Definition 2.1 the \emph{log structure} on $\Spec(R)$ attached to $P\ra R$
is \emph{log regular}. Here, this means that the residue class ring $\bar{R}=R/(P\smallsetminus P^\times)$ is regular,
which is a consequence of Lemma \ref{maximal ideal}, and 
$$
\dim(R)=\dim(\bar{R})+\rank_\ZZ(P^\gp),
$$
which holds because the numbers are $2=0+2$. We refrain from recalling more details, since
the full strength of Kato's theory is not needed, but   stress again that his ideas
play a paramount role throughout.  From loc.\ cit., Theorem 3.1 and Theorem 3.2, it follows:

\begin{proposition}
\mylabel{kernel of surjection}
The kernel of the surjection  $W[[P]]\ra R$ is
generated by a formal power series
of the form
$$
\psi=\begin{cases}
p + \sum_{g\in P\smallsetminus P^\times} w_g  g & \text{if $R$ does not contain a field;}\\
0                                            & \text{if $R$ contains a field.}
\end{cases}
$$
In particular, the map $W[[P]]\ra R$ is bijective if the ring $R$ contains a field.
\end{proposition}

Consider the commutative diagram of monoids and groups
$$
\begin{CD}
P @>>> Q @>>>\NN\oplus\NN\\
@VVV  @VVV @VVV\\
P^\gp @>>> Q^\gp @>>> \ZZ\oplus\ZZ
\end{CD}
$$
Here $P^\gp\subset Q^\gp\subset\ZZ\oplus\ZZ$ is the unique intermediate group whose index $[\ZZ\oplus\ZZ:Q^\gp]$
is a $p$-power, and that the other index $[Q^\gp:P^\gp]$ is prime to $p$.
Recall once more that $p\geq 1$ denotes the  characteristic exponent of the residue field $R/\maxid_R$.
The monoid $Q$ is defined by demanding that the square on the right is cartesian.
This $Q$ is sharp and fs, compare the arguments for Corollary \ref{properties P}. We now define rings
\begin{gather*}
R':=R\otimes_{W[[P]]}W[[Q]]=W[[Q]]/(\psi),\\
R'':=R\otimes_{W[[P]]}W[[\NN\oplus\NN]]=W[[\NN\oplus\NN]]/(\psi).
\end{gather*}
This gives finite   extensions $R\subset R'\subset R''$ between complete local noetherian rings.

\begin{proposition}
\mylabel{ring properties}
Both  rings $R',R''$ are   integral domains. The ring $R'$ is normal,   the ring $R''$ is regular,
the ring extension $R\subset R'$ is \'etale in codimension one, and the ring extension $R'\subset R''$ is totally ramified
at some prime $\mathfrak{q}\subset R'$ of height one.
\end{proposition}

\proof
The rings in question are integral by \cite{Kato 1994}, Lemma 3.4.
Endowed with their canonical log structures, they become log regular, whence are normal,
according to loc.\ cit., Theorem 3.1 and Theorem 4.1.
Obviously, $W[[\NN\oplus\NN]]$ is regular.
Its residue class ring $R''$ must be regular as well, in light of the description in Proposition \ref{kernel of surjection}.
 
Set 
$$
Y=\Spec(W[[P]])\quadand Y'=\Spec(W[[Q]])\quadand Y''=\Spec W[[\NN\oplus\NN]],
$$
and consider the ensuing  morphisms $Y''\ra Y'\ra Y$.
The ring $W[[\NN\oplus\NN]]$ is endowed with a natural grading by the finite cyclic group
$N=(\ZZ\oplus\ZZ)/P^\gp$: The formal power series $\sum w_gg$, with $g\cong d$ modulo $P^\gp$
are the homogeneous elements of degree $[d]\in N$. In turn, we get an action
of the finite flat diagonalizable group scheme $H=\Spec W[N]$ on $Y''$ with
quotient $Y=Y''/H$, as explained in \cite{SGA 3b}, Expose VIII.
According to Proposition \ref{free action}, the quotient morphism $Y''\ra Y$ 
is a $H$-torsor over the complement of the closed subscheme $\Spec(W)\subset\Spec(W[[P]])$
defined by the ideal $(P\smallsetminus P^\gp)\subset W[[P]]$.

Now consider the finite flat subgroup scheme $H'\subset H$ corresponding
to the finite cyclic group $N'=(\ZZ\oplus\ZZ)/Q^\gp$, such that $Y'=Y''/H'$.
Restricted to the residue field $W/\maxid_W$, the group scheme $H'=\Spec W[N']$
becomes radical, because   this residue field has characteristic exponent $p$
and the group $N$ has order a power of $p$.
We conclude that $R'\subset R''$ must be totally  ramified along some prime of height one.
Note that this holds for trivial reasons if  $R$ contains a field of characteristic zero, for then $R'=R''$.

Finally, consider the finite flat quotient group scheme $G=H/H'$ corresponding to the
finite cyclic group $M=Q^\gp/P^\gp$, such that $Y=Y'/G$.
The group scheme $G$ is \'etale, because the order of $M$ is prime to
the characteristic exponent of the residue field $W/\maxid_W$.
Using Proposition \ref{free action} again, we easily deduce that the quotient morphism $Y'\ra Y$
is a $G$-torsor over the complement of $\Spec(W)\subset\Spec W[[P]]$
defined by the ideal $(P\smallsetminus P^\gp)\subset W[[P]]$.
From the description of the closed subscheme $\Spec(R)\subset\Spec W[[P]]$ in Proposition \ref{kernel of surjection},
we infer that $R\subset R'$ must be a $G$-torsor in codimension one, and in particular
is \'etale in codimension one.
\qed

\medskip
Let $R\subset K$ and $R'\subset K'$ be the field of fractions.

\begin{proposition}
\mylabel{field properties}
The finite field extension $K\subset K'$ is cyclic, and $\Gal(K'/K)$
is isomorphic to the   prime-to-$p$ part of $\ZZ/m\ZZ$.
\end{proposition}

\proof
Let $\ZZ/m'\ZZ$ be the   prime-to-$p$ part of $\ZZ/m\ZZ$.
This group is isomorphic to $M=Q^\gp/P^\gp$, by definition of $Q^\gp$.
Since $m'$ is invertible in $W$ and the residue field $W/\maxid_W$ is separably closed, 
the group scheme $G$ representing $A\mapsto\Hom(M,A^\times)$
comes from the abstract group $\ZZ/m'\ZZ$. We say in the proof of the previous proposition
that $\Spec(K')\ra\Spec(K)$ is a $G$-torsor, in other words, $K\subset K'$ is Galois,
with Galois group $\Gal(K'/K)\simeq\ZZ/m'\ZZ$.
\qed

\medskip
\emph{Proof of Theorem \ref{fundamental group}.}
Using Proposition \ref{galois group}, together with Proposition \ref{ring properties} and \ref{field properties},
we find that $\pi_1^\loc(R)=\Gal(K'/K)$ is isomorphic to the   prime-to-$p$ part
of $\ZZ/m\ZZ$, where $R$ is a Hirzebruch--Jung singularity of type $m/b$.
\qed

%===========================================================
\section{Cotangent spaces}
\mylabel{Cotangent spaces} 

The goal of this section is to prove Lemma \ref{maximal ideal}, which says that the
maximal ideal of   Hirzebruch--Jung singularities is generated by
the nonunits of certain monoids.

We start with some preparatory considerations. Let $k$ be a field, and $E$ be 
a proper reduced curve over $k$ whose irreducible components $E_1,\ldots, E_r$
are  isomorphic to $\PP^1_k$, with 
$$
h^0(\O_{E_i\cap E_j}) =\dim_k H^0(\O_{E_i\cap E_j})= \begin{cases}
1 & \text{if $|i-j|=1$;}\\
0 & \text{if $|i-j|\geq 2$.}
\end{cases}
$$
To rule out triple intersections, we  furthermore demand that each point $x\in E$ lies in at most two irreducible components.
We also consider    noetherian reduced 1-dimensional schemes having the form
$$
\tilde{E} = E\cup D\cup D',
$$
where $D,D'$ are the spectra of certain discrete valuation rings, not necessarily containing $k$,
so that $h^0(\O_{E_i\cap D})=\delta_{i,0}$ and $h^0(\O_{E_i\cap D'})=\delta_{i,r}$ (Kronecker delta).
As we saw in the preceding section, such curves occur  on the resolution of singularities for  Hirzebruch--Jung
singularities. To make formulas more elegant, we also write $E_0=D$ and $E_{r+1}=D'$.
Thus the dual graph of $\tilde{E}$ is:

%-------------------------------------------------------------
\vspace{2em}
\centerline{\includegraphics{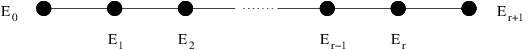} }
\vspace{1em}
%\centerline{Figure \stepcounter{figure}\arabic{figure}: The dual graph of $\tilde{E}$.}
\vspace{1em}
%-------------------------------------------------------------
Now let $\shL$ be an invertible sheaf on $E$ with intersection numbers
$$
d_i=(\shL\cdot E_i) = \begin{cases}
\geq 0 & \text{if $i=1$ or $i=r$;}\\
\geq 1 & \text{if $2\leq i\leq r-1$}.
\end{cases}
$$
Here the intersection numbers are defined as $(\shL\cdot E_i)=\chi(\shL_{E_i})-\chi(\O_{E_i})$ 
via $k$-vector space dimensions. For the proof of Lemma \ref{maximal ideal}, 
it will be useful to have a specific basis for the finite dimensional $k$-vector space $H^0(E,\shL)$
adapted to our problem. The dimension is easy:

\begin{proposition}
\mylabel{global sections}
We have $h^0(\shL)=1+\sum_{i=1}^r d_i$.
\end{proposition}

\proof

The dualizing sheaf $\omega_E$ is invertible, since $E$ is locally of complete intersection
(by \cite{EGA IVd}, Corollary 19.3.4, it suffices to check this over the algebraic closure of $k$,  where
it is well-known), and has intersection numbers
$$
(\omega_E\cdot E_i) = \begin{cases}
-1 & \text{if $i=1$ or $i=r$;}\\
0  & \text{if $2\leq i\leq r-1$.}
\end{cases}
$$
Thus $(\shL^\vee\otimes\omega_E \cdot E_i)<0$. By Serre duality, $H^0(E,\shL^\vee\otimes\omega_E)\simeq H^1(E,\shL)$,
and the latter vanishes since $E$ is reduced. Whence
$$
h^0(\shL)=\chi(\shL) =\deg(\shL) +\chi(\O_E) = (\sum_{i=1}^r d_i) +1
$$
from Riemann--Roch.
\qed

\medskip
We now define $k$-rational points
$$
a_i,b_i\in E_i,\qquad 1\leq i\leq r
$$
by the condition $E_i\cap E_{i-1} =\left\{a_i\right\}$ and $E_i\cap E_{i+1}=\left\{b_i\right\}$.
Note that we have $a_i=b_{i+1}$ as points in $E=E_1\cup\ldots \cup E_r$. We also set  
$b_0=a_1$ and $a_{r+1}=b_r$, and call the  rational points
$$
a_1=b_0,\quad a_2=b_1,\quad \ldots \quad  a_{r+1}=b_r\in E
$$
the \emph{special points} of the curve $E$. These are precisely the singular points on  the extended curve $\tilde{E}$. 
Using these special points, one obtains an exact sequence
\begin{equation}
\label{section sequence}
0\lra H^0(E,\shL) \lra \bigoplus H^0(E_i,\shL|_{E_i}) \lra \bigoplus_{i=1}^{r-1} \frac{\kappa(b_i)\oplus\kappa(a_{i+1})}{k},
\end{equation}
where the map on the left is the sum of the restriction maps,   the quotients on the right are with
respect to the diagonal inclusions of the ground field in the residue fields, and the map on the right
is induced from the evaluation  map
$$
(f_1,\ldots,f_r)\longmapsto (f_1(b_1),f_2(a_2),f_2(b_2),f_3(a_3),\ldots,f_{r-1}(b_{r-1}),f_r(a_r)).
$$
Here $f_i\in H^0(E_i,\shL|_{E_i})$. Now for each $1\leq i\leq r$, choose $d_i-1$ sections
$$
s_{ij}\in H^0(E_i,\shL|_{E_i}), \qquad 1\leq j\leq d_i-1
$$
so that $s_{ij}$ vanishes at the special point $a_i\in E_i$ of order $j$, and at the special
point $b_i\in E_i$ of order $d_i-j$. For example, if one writes $E_i=\Proj k[T_0,T_1]$ as a homogeneous spectrum
and identifies $\shL|_{E_i}$ with $\O_{\PP^1}(d_i)$,
then one may choose $s_{ij}=T_0^jT_1^{d_i-j}$.
Using the exact sequence (\ref{section sequence}), we may  regard these sections over $E_i$ also as global sections
$s_{ij}\in H^0(E,\shL)$, because the $s_{ij}$ vanish at all special points.

\begin{proposition}
\mylabel{basis special}
The $\sum_{i=1}^r(d_i-1)$ sections
$$
s_{ij}\in H^0(E,\shL),\quad 1\leq i\leq r,\quad 1\leq j\leq d_i-1
$$
form a basis for the vector subspace $V\subset H^0(E,\shL)$ of all sections vanishing
at the special points.
\end{proposition}

\proof
Using that the locus on nonvanishing $E_{s_{ij}}=E_{i}\smallsetminus\left\{a_i,b_i\right\}$ 
are pairwise disjoint for $1\leq i\leq r$,
and that the $s_{ij}\in H^0(E_i,\shL|_{E_i})$ are linearly independent for $1\leq j\leq d_i-1$, one infers
that the $s_{ij}\in H^0(E,\shL)$ are linearly independent. On the other hand, we have  an inclusion
$$
V\subset\bigoplus H^0(E_i,\shL|_{E_i}(-a_i-b_i)),
$$
and the right hand side is a vector space of dimension $\sum_{i=1}^r(d_i-1)$. If follows that the $s_{ij}\in V$ form a basis.
\qed

\medskip
Next, suppose we have sections 
$$
t_i\in H^0(E,\shL),\quad 1\leq i\leq r+1
$$
with $t_i(a_i)\neq 0$ in the residue field $\kappa(a_i)$, and $t_i|_{E_n}=0$ as elements of $H^0(E_i,\shL|_{E_i})$
for all indices $n\neq i,i-1$. This will give us the desired basis:

\begin{proposition}
\mylabel{basis}
Assumptions as above. The set
$$
\left\{s_{ij}\mid 1\leq i\leq r,\, 1\leq j\leq d_i-1\right\}\quad\cup\quad
\left\{t_i\mid 1\leq i\leq r+1\right\}   
$$
comprise a basis for the $k$-vector space $H^0(E,\shL)$.
\end{proposition}

\proof
The set has cardinality $\sum_{i=1}^r(d_i-1) +(r+1)=h^0(\shL)$, so it suffices to check that the sections
form a generating system. Let $s\in H^0(E,\shL)$ be arbitrary. Subtracting a suitable linear
combination of the $t_i$, we may assume that $s$ vanishes at the special points.
Then our $s$ becomes a linear combination of the $s_{ij}$ by Proposition \ref{basis special}.
\qed

\medskip
Now let $R$ be a Hirzebruch--Jung singularity, with minimal resolution of singularities
$f:X\ra\Spec(R)$.  
The reduced exceptional divisor $E\subset X$ is a proper reduced curve over the
residue field $k=R/\maxid_R$ as studied above.
In the preceding section, we have chosen two integral Cartier divisors $D,D'\subset X$,
which give us the   1-dimensional scheme $\tilde{E}=E\cup D\cup D'$, and whose images are
integral Weil divisors $C,C'\subset\Spec(R)$. From this, we have constructed
a sharp fs submonoid $P\subset R$. Our task now is to prove Lemma \ref{maximal ideal},
that is $P\smallsetminus P^\times$ generates the maximal ideal $\maxid_R$.

To do so we recall Artin's description of the cotangent space $\maxid_R/\maxid_R^2$
in term of the fundamental cycle $Z\subset X$ of a rational singularity in dimension two. 
The definition of the fundamental cycle ensures
that the ideal $\maxid_R\O_X\subset\O_X$
is contained in the ideal $\O_X(-Z)\subset\O_X$. Using the short exact sequence
$$0\ra\O_X(-2Z)\ra\O_X(-Z)\ra\O_Z(-Z)\ra 0,$$ one   obtains a   map
$$
\maxid_R/\maxid_R^2\lra H^0(Z,\O_Z(-Z)).
$$
It follows from \cite{Artin 1966}, Theorem 4 that this map is bijective.
In our Hirzebruch--Jung situation, the fundamental cycle is $Z=E$. Setting $\shL=\O_E(-E)$, we have degrees
$$
d_i=(\shL\cdot E_i) = \begin{cases}
s_i-2 & \text{for $2\leq i\leq r-1$;}\\
s_i-1 & \text{for $i=1$ or $i=r$},
\end{cases}
$$
with the  numbers $s_i=-E_i^2\geq 2$.
Now consider an effective divisor $A=\sum_{n=0}^{r+1}\lambda_nE_n$ that is numerically
trivial on $E$. The latter is  equivalent to the system of equations
\begin{equation}
\label{trivial intersections}
s_n\lambda_n = \lambda_{n-1}+\lambda_{n+1},\quad 1\leq n\leq r.
\end{equation}
In particular, we have the mean value estimates
\begin{equation}
\label{mean values}
\lambda_n \leq \frac{\lambda_{n-1}+\lambda_{n+1}}{2},\quad 1\leq n\leq r.
\end{equation}
Since $R$ is rational, the ideal of the closed subscheme $A\subset X$ is of the form $\O_X(-A)=g\O_X$ for some
$g\in \maxid_R$, which is unique up to $R^\times$,
and in turn gives a global section $\bar{g}\in H^0(E,\shL)$, the latter being unique up to $k^\times$.
Now recall the construction of our basis $s_{ij},t_i\in H^0(E,\shL)$ from Proposition \ref{basis}.

\begin{proposition}
\mylabel{first sections}
Notation as above. Suppose for some $1\leq i\leq r$, we have $\lambda_i=1$ and $\lambda_n\geq 2$
for $n=1,\ldots,i-1$ and $n=i+1,\ldots,r$.
Then $\bar{g}\in H^0(E,\shL)$ is a nonzero multiple of $s_{ij}$, with $j=\lambda_{i-1}-1$.
\end{proposition}

\proof
It is clear that $\bar{g}$ vanishes on each $E_n$, $n\neq i$, but does not vanish
on $E_i\smallsetminus\left\{a_i,b_i\right\}$. It thus suffices to check that $\bar{g}$
vanishes of order $j$ at $a_i\in E_i$. This is a local problem: Set $x=a_i$, and let $u,v\in \O_{X,x}$
be a   system of parameters whose zero locus is $E_{i-1},E_i$, respectively.
Up to units, we have $g=u^j\cdot uv$.
Under the canonical surjection $\O_X(-E)\ra\O_{E_i}(-E)$, the element $g$ becomes the class of $u^j$
in $\O_{X,x}/(v)=\O_{E_i,x}$.
\qed

\medskip
In the same way, one checks:

\begin{proposition}
\mylabel{second sections}
Notation as above. Suppose for some index $1\leq i\leq r+1$ we have
$\lambda_{i-1}=\lambda_{i}=1$ and
$\lambda_n\geq 2$ for $n=1,\ldots, i-1$ and $n=i+1,\ldots,r$.
Then $\bar{g}\in H^0(E,\shL)$ is a nonzero multiple of $t_i$.
\end{proposition}

We   come to the final goal of this section:

\medskip\noindent\emph{Proof of Lemma \ref{maximal ideal}.} 
By the Nakayama Lemma, it suffices to check that the image of $P\smallsetminus P^\times$
generates the vector space $\maxid_R/\maxid_R^2=H^0(E,\shL)$, for the invertible sheaf $\shL=\O_E(-E)$.

Fix a basis element $s_{ij}\in H^0(E,\shL)$,   for some $1\leq i\leq r$ and $1\leq j\leq d_i-1$.
We   define an effective Cartier divisor $A=\sum_{i=0}^{r+1}\lambda_iE_i$ as follows:
First set
$$\lambda_i=1\quadand\lambda_{i-1}=j+1\quadand\lambda_{i+1}=s_i-1-j,
$$
and  then define the preceding multiplicities $\lambda_{i-2},\ldots,\lambda_0$ by descending
induction through the equations (\ref{trivial intersections}),
and similarly for the subsequent $\lambda_{i+2},\ldots,\lambda_{r+1}$.
The mean value estimates (\ref{mean values}) ensure that the sequence of integers $\lambda_0,\ldots,\lambda_{r+1}$ is
strictly decreasing until its minimum $\lambda_i=1$, after which it becomes strictly increasing.
By construction, the effective Cartier divisor $A=\sum_{i=0}^{r+1}\lambda_iE_i$ is numerically trivial on $E$ and  has
multiplicities $\geq 1$,  thus yields an element $g\in P\smallsetminus P^\times\subset\maxid_R$.
In light of Proposition \ref{first sections}, the induced section $\bar{g}\in H^0(E,\shL)$ is a nonzero scalar multiple of $s_{ij}$.

Now fix a basis element $t_i\in H^0(E,\shL)$. Then we define an effective Cartier divisor $A=\sum_{i=0}^{r+1}\lambda_iE_i$
in a similar way: First set
$$
\lambda_{i-1}=\lambda_i=1\quadand
\lambda_{i-2}=s_{i-1}-1\quadand\lambda_{i+1}=s_i-1,
$$
and define the remaining multiplicities as
in the preceding paragraph.  We infer with Proposition \ref{second sections}   that $\bar{g}$ is a nonzero scalar multiple of $t_i$.

Proposition \ref{basis} tells us that the image of  $P\smallsetminus P^\times\ra H^0(E,\shL)$ 
contains a basis. Thus  $P\smallsetminus P^\times$ generates the maximal ideal $\maxid_R$.
\qed

%===========================================================
\section{Hermitian curves and special $p$-groups}
\mylabel{Hermitian curves} 

In the remaining sections, we investigate certain curves with very large automorphism groups,
in order to produce new examples of wild quotient surface singularities whose dual graph contains
at least two nodes.

Let $p>0$ be a prime number, $k$ an algebraically closed ground field of characteristic $p$, and $q=p^m$ be a fixed prime
power. Consider the smooth projective curve given by the affine equation
$$
C:\quad y^q -y = x^{q+1},
$$
which has genus $g=q(q-1)/2$. These are the so-called \emph{Hermitian curves}, and are known for extremal behavior
with respect to automorphisms and rational points, compare   \cite{Stichtenoth 1973a}, \cite{Rueck; Stichtenoth 1994}.
Shioda kindly informed us that this curve is isomorphic to the Fermat curve of degree $q+1$, which are discussed
in \cite{Shioda 1988}, an isomorphism being defined over the field $\FF_{q^2}$. We prefer, however, to work with 
the equation for the Hermitian curve.

The substitution $x = x'+r$ and  $y = y'+sx'+t$ transforms the equation into
$$
y'^q-y' = x'^{q+1} +(r-s^p)x'^q + (r^q+s)x' + (r^{q+1}-t^q+t).
$$
Such a substitution leaves the original equation of the Hermitian curve invariant if and only if
$r-s^q=0$, $r^q+s=0$, $r^{q+1}-t^q+t=0$, which is equivalent to 
$$
r^{q^2}+r=0\quadand t^q-t = r^{q+1}\quadand s=-r^q.
$$
Summing up, we may regard the set
$$
G=\left\{(t,r)\in k^2\mid \text{$r^{q^2}+r=0$ and $t^q-t = r^{q+1}$}\right\}
$$
as a group of automorphisms, where the elements act  via
\begin{equation}
\label{group action}
(t,r):\quad x\longmapsto x + r \qquad y\longmapsto y -r^qx + t.
\end{equation}
The group law, viewed as composition of functions, is 
$$
(t,r) \circ (t',r') = (t+t'-rr'^q,r+r').
$$
However, we endow the set $G$ with the opposite group law
$$
(t,r) \cdot (t',r') = (t+t'-r^qr',r+r')
$$
since we want to regard $G$ as a group of automorphisms of the scheme $C$ rather than 
its   coordinate ring.

% 
% It is clear from this description that a nonzero group element does not fix any point
% on the affine part of $C$. Thus:
% 
% \begin{proposition}
% \mylabel{free action}
% The $G$ action on the affine part of $C$ is free.
% \end{proposition}

\begin{proposition}
\mylabel{order}
The  group $G$ has order $\ord(G)=q^3$. It is a Sylow-$p$-subgroup of $\Aut(C)$ when    $q\neq 2$.
\end{proposition}

\proof
Polynomials over $k$ of the form $T^{q^2}+T$ or $T^q-T-\lambda$ are separable, so the statement about the order follows.
For $q\neq 2$, the full automorphism group $\Aut(C)$ has order $q^3(q^3+1)(q^2-1)$, according to \cite{Stichtenoth 1973a}, Hauptsatz. 
\qed

\medskip
Note that for $q=2$,
the curve $C$ is the supersingular elliptic curve in characteristic two, and then the group of automorphisms fixing the origin
has order $24$. So  $G$ is a Sylow-$p$-subgroup as well.

Recall that the \emph{Frattini subgroup} of a group is the intersection of all maximal proper subgroups.
For our given group $G$, we denote by $\Phi,G',Z\subset G$   the Frattini subgroup, the commutator subgroup,
and the center, respectively.

\begin{proposition}
\mylabel{special p-group}
In our group $G$,   Frattini subgroup,   commutator subgroup and center coincide: $\Phi= G'=Z$.
This subgroup is given by $\left\{(t,0)\mid t\in\FF_q\right\}\subset G$.  
\end{proposition}

\proof
A direct computations shows that the centralizer of an element $(t,r)\in G$ with $r\neq 0$ consists of
those $(t',r')$ with $r'=\lambda r$ for some $\lambda\in\FF_q$. It follows that
$Z=\left\{(t,0)\mid t\in\FF_q\right\}$, and this is obviously the kernel of the homomorphism
$$
G\lra k,\quad (t,r)\longmapsto r.
$$
The inclusion $G'\subset\Phi$ holds
for arbitrary finite $p$-groups, because then the Frattini subgroup is generated by the commutators and
$p$-powers (see   \cite{Macdonald 1968}, Theorem 9.26). Thus  $G/\Phi$ is the largest elementary
abelian quotient, and we conclude $\Phi\subset Z$.

It remains to check that each central element $(t,0)$ is a commutator. Indeed,
inverses and commutators   are given by the respective formulas
\begin{equation}
\label{inverse and commutator}
(t,r)^{-1} = (-r^{q+1}-t,-r)\quadand [(t,r),(t',r')]=(r^qr'-rr'^q,0).
\end{equation}
Now observe that the formula for  commutators is $\FF_q$-linear in $r$ and $r'$  with respect to the   $\FF_q$-vector space structures. 
Since there is at least one nontrivial commutator, each element $(t,0)$ is automatically a commutator.
\qed

\medskip
In other words, our group $G$ is  a \emph{special $p$-group}, which by definition means that   $\Phi=G'=Z$.
These groups play a role in the classification of finite simple groups. It   follows from the definition 
that abelianization $G^\ab=G/G'$ is    elementary abelian,
and the same holds for the center $Z\subset G$, see \cite{Aschbacher 1986}, (23.7).
Moreover,  the pairing
$$
G^\ab\times G^\ab\lra Z,\quad (\bar{a}, \bar{b})\longmapsto [a,b]=aba^{-1}b^{-1}
$$
is  nondegenerate and alternating. 

In the extremal case where $\Phi=G'=Z$ is cyclic, hence of order $p$, the group $G$ is called
\emph{extraspecial}; then the pairing is a symplectic form, and one has $\ord(G)=p^{1+2n}$ for some $n\geq 0$.
There is a classification of extraspecial $p$-groups: For each order, there are precisely two isomorphism classes,
compare \cite{Doerk; Hawkes 1992}, Chapter B, \S9. 
For odd $p$, they are distinguished by the exponent of the group.
The situation is somewhat more complicated for $p=2$: The two extraspecial $2$-groups of order eight are
the dihedral group $D_4=C_4\rtimes\left\{\pm 1\right\}$ and the quaternion group $Q=\left\{\pm 1,\pm i,\pm j,\pm k\right\}$.
In our situation, we have:

\begin{proposition}
\mylabel{exponent}
Our special $p$-group $G\subset\Aut(C)$ has exponent
$$
\exp(G)=
\begin{cases}
p   & \text{if $p\neq 2$};\\
p^2 & \text{if $p=2$}. 
\end{cases}
$$
Moreover, for $p=2$ all noncentral elements in $G$ have order four.
\end{proposition}

\proof
Since $Z$ and $G^\ab$ are elementary abelian,  the exponent divides $p^2$. One easily computes that $p$-th powers are given by
$(t,r)^p= (\sum_{i=1}^{p-1}i r^{p+1},0)$.
The first entry is $\sum_{i=1}^{p-1}i r^{p+1}= r^{p+1} \cdot p(p-1)/2$, and the result follows.
\qed

%===========================================================
\section{Quotient singularities with two nodes}
\mylabel{Two nodes}

Fix a prime $p>0 $ and let $q=p^m$ be a prime power. Let  $C:\;y^q-y=x^{q+1}$ be the Hermitian curve
and $G$ be the special $p$-group  of order $\ord(G)=q^3$ acting on $C$, as described in the previous section.
The $G$-action is free on the affine part of $C$, and the point $\infty\in C$ at infinity
is fixed.
We now consider the diagonal action of $G$  on the smooth proper surface $C\times C$.
Form the quotient surface
$$
Y=G\backslash(C\times C).
$$
This is a normal proper surface,
whose singular locus is the point $y\in Y$ corresponding to the fixed point $(\infty,\infty)\in C\times C$. 
Our second main result is a statement on the dual graph of this wild quotient singularity,
which gives a partial answer to Lorenzinis Question 1.1 b) in \cite{Lorenzini 2011d}:

\begin{theorem}
\mylabel{two nodes}
The dual graph for the minimal resolution of the wild quotient singularity $y\in Y=G\backslash(C\times C)$
contains at least two nodes.
\end{theorem}

Note that we have no explicit description of this wild quotient singularity in terms of formal equations, nor do we know the
precise number of nodes or the exact shape of the dual graph.
We shall deduce the theorem by comparing with another singularity that can be handled.
To this end we consider the fibration
$$
\varphi:Y=G\backslash(C\times C)\lra G\backslash C=\PP^1
$$
induced from the   projection onto the first factor. This fibration admits
a section, given by  $ \left\{\infty\right\}\times C$. Now consider the commutative diagram 
$$
\begin{xy}
\xymatrix{
                & \tilde{Y}\ar[dl]\ar[dr] \\
S\ar[dr]_{\psi} &                         & Y=G\backslash(C\times C)\ar[dl]^{\varphi} & C\times C \ar[l]\ar[d]^{\pr_1} \\
                & \PP^1=G\backslash C     &                                           & C,\ar[ll]
}
\end{xy}
$$
where $\tilde{Y}\ra Y$ is the minimal resolution of singularities, and $\psi:S\ra\PP^1$ is the relative minimal model,
that is, the model obtained by successively contracting $(-1)$-curves in the singular fiber of $\tilde{Y}\ra\PP^1$.
By abuse of notation, we write $\infty\in\PP^1$ for
the point corresponding to the fixed point $\infty\in C$, such that $\psi^{-1}(\infty)\subset S$ is   the singular fiber.

%-------------------------------------------------------------
\vspace{2em}
\centerline{\includegraphics{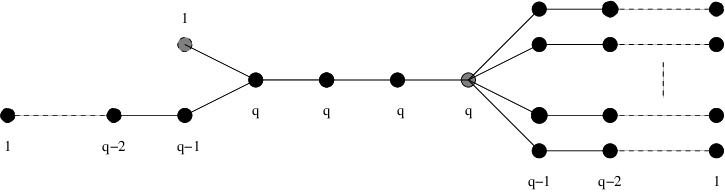}} 
\vspace{1em}
\centerline{Figure \stepcounter{figure}\arabic{figure}: Dual graph of $\psi^{-1}(\infty)$ with multiplicities}
\vspace{1em}
%-------------------------------------------------------------

\begin{proposition}
\mylabel{singular fiber}
The reduced singular fiber  $\psi^{-1}(\infty)_\red\subset S$ is a divisor with strictly normal crossings whose irreducible components
are copies of $\PP^1$, and its dual graph
is depicted in Figure 1. 
There are $q$ strings on the right hand side, each of length $q-1$.
The numbers indicate   multiplicities of integral components in the schematic fiber $\psi^{-1}(\infty)$.
Integral components corresponding to the black vertices have self-intersection $-2$, whereas the other two
have self-intersection number $-q$.
\end{proposition}

For this we proceed in two steps: First, we determine an explicit equation for the generic fiber
$\psi^{-1}(\eta)$, where $\eta\in\PP^1=G\backslash C$ denotes the generic point. Second, we use this equation to find an explicit resolution  
for the ensuing singularities, which finally will produces the dual graph $\psi^{-1}(\infty)$. The arguments
will occupy the next two sections.

\medskip\noindent
\emph{Proof that Proposition \ref{singular fiber} implies Theorem \ref{two nodes}.}
Consider the induced fibration $\tilde{Y}\ra\PP^1$. Its reduced singular fiber consists of two parts:
On the one hand, the strict transform $F\subset \tilde{Y}$ of the singular fiber of $\varphi:Y\ra\PP^1$.
On the other hand, the exceptional divisor $E\subset\tilde{Y}$ coming from the resolution of singularities.

The multiplicity of the integral component $F$ in the schematic fiber is $q^3=\ord(G)$, 
compare \cite{Ito; Schroeer 2012}, Proposition 2.1.
By Proposition \ref{singular fiber}, there are no integral components of multiplicity $>q$ in the 
singular fiber $\psi^{-1}(\infty)\subset S$. Whence $F$ gets contracted by $\tilde{Y}\ra S$.
Therefore, the dual graph of $\psi^{-1}(\infty)$ is obtained from the dual graph of $E$ by
a sequence of vertex contractions.
Again by Proposition \ref{singular fiber}, the dual graph for $\psi^{-1}(\infty)_\red$ contains
precisely two nodes. This obviously  implies that the dual graph for $E$ contains at least two nodes.
\qed

%===========================================================
\section{Description of the generic fiber}
\mylabel{Generic fiber}

Our goal here is to describe the generic fiber $\varphi^{-1}(\eta)$
of the fibration 
$$
\varphi: G\backslash(C\times C)\ra G\backslash C=\PP^1
$$
induced from projection
onto the first factor. Recall that $q=p^m$ is our fixed prime power, and
$C:\;y^q-y=x^{q+1}$ is our Hermitian curve on which the group $G$ of order $\ord(G)=q^3$ acts.

Let $R$ be an arbitrary ring of characteristic $p$, and consider the $R$-algebra
$$
A=R[X,Y]/(Y^q-Y-X^{q+1}),
$$
endowed with the canonical $G$-action as in  (\ref{group action}). 
Here we temporarily use upper case symbols $X,Y$ in order to distinguish later
the two factors  in the selfproduct $C\times C$. Let   $\tau\in A$ be the residue class of $X^{q^2}+X$.
From the description of the  action  one sees that $\tau$ is  $G$-invariant.

\begin{proposition}
\mylabel{G-invariants}
The subring $R[\tau]\subset A$ is the ring of $G$-invariants.
\end{proposition}

\proof
Considering the intermediate ring $R[X]$, one sees that
the ring   $A$ is   free of rank $q^3=\ord(G)$ as module over $R[\tau]$.
To finish the argument, it suffices to treat the case that $R$ is noetherian. By the Nakayama Lemma,
we are reduced to the case that $R$ is a field. Then both $R[\tau]$ and $A$ are integral and normal.
Now Galois theory and Zariski's Main Theorem yield the result.
\qed

\medskip
We  view $U=\Spec k[\tau]$ as the open subset of $\PP^1=G\backslash C$ over which the fibration
$\varphi:G\backslash(C\times C)\ra\PP^1$ is smooth. In particular, the   preimage $\varphi^{-1}(U)$ becomes
a flat family of smooth projective curves over $U$.  

\begin{proposition}
\mylabel{twisted form}
The   family of smooth projective curves $\varphi^{-1}(U)\ra U$  is given by the homogenization of the equation 
$y^q-y=x^{q+1}+\tau x^q$.
\end{proposition}

\proof
Let us point out that the challenge was to \emph{guess} the right equation.
Once an equation is found, it is  
not too difficult to verify its correctness with the abstract theory of 
twisted forms, as exposed at length by Giraud in his monograph \cite{Giraud 1971}.

Consider the smooth projective families of curves $\shC\ra U$ and $\shC'\ra U$ given by the  homogenizations of the
equations
$$
y^q-y=x^{q+1}  \quadand y'^q-y'=x'^q+\tau x'^q,
$$
respectively. To compare them we use the  functor $\underline{\Isom}(\shC,\shC')$, whose values  on  $k[\tau]$-algebras $R$ 
is the set of $R$-isomorphisms $\shC\otimes_{k[\tau]} R\ra\shC'\otimes_{k[\tau]} R$. It contains as a subfunctor 
$$
I(R)= \left\{(b,a)\in R^2\mid a^{q^2}+a+\tau=0\quadand b-b^q+\tau^{q+1}-a^{q^2+q}=0\right\},
$$
where the inclusion $I(R)\subset\underline{\Isom}(\shC,\shC')(R)$ comes from the substitutions $x'\mapsto x+a$ and $y'\mapsto y-a^qx+b$. 
To proceed, we consider another functor
$$
P(R)=\left\{(b,a)\in R^2\mid a^{q^2}+a-\tau = 0\quadand b^q-b-a^{q+1}=0\right\},
$$
which is representable by the spectrum of $k[\tau,X,Y]/(X^{q^2}+X-\tau,Y^q-Y-X^{q+1})$.
Note that the latter is a finite \'etale $k[\tau]$-algebra. We have a natural transformation
$$
\Upsilon:P(R)\lra I(R), \quad (b,a)\longmapsto(-b-a^{q+1},-a).
$$
Our   special $p$-group
$$
G=\left\{(t,r)\in k^2\mid r^{q^2}+r=0\quadand t^q-t=r^{q+1}\right\}\subset\Aut(C)
$$
acts on the right on $\underline{\Isom}(\shC,\shC')$ by composition, and one easily computes that
this action is given by
$$
(b,a)\circ(t,r) = (b+t-a^qr,a+r).
$$
On the other hand, $P$ is endowed with a canonical $G$-action from the left, given by
$(t,r)\cdot(b,a)=(b-r^qa+t,a+r)$ as in (\ref{group action}). In fact, the structure morphism
$P\ra U$ becomes a   $G$-torsor. The corresponding \emph{right} action is
$$
(b,a)\cdot(t,r) = (t,r)^{-1}\cdot(b,a) = (-r^{q+1}-t+b+r^qa,a-r),
$$
compare (\ref{inverse and commutator}). A straightforward computation reveals that the natural transformation
$\Upsilon:P\ra I$ respects the right $G$-actions. Actually, it is a $u$-morphism, where
$u:G\ra\underline{\Aut}(\shC)$ denotes the canonical inclusion. In turn, we obtain a morphism 
$$
P\wedge^G\underline{\shC}\lra \underline{\Isom}(\shC,\shC')
$$
of right $\underline{\Aut}(\shC)$-torsors. Being a morphism of torsors, it is automatically an isomorphism.
According to   \cite{Giraud 1971}, Chapter III, Theorem 2.5.1,
the family of curves $\shC'$ is a twisted form of $\shC$ with respect to the \'etale topology on $U$,
and one has $\shC'\simeq P\wedge^G\shC$. The latter denotes the \emph{contracted product} and is defined as
 the quotient of $P\times \shC$ by the left $G$-actions
 $g\cdot (p,c)=(pg^{-1},gc) = (gp,gc)$.
But this coincides with the very definition of $\varphi^{-1}(U)\subset G\backslash (C\times C)$.
The upshot is the $\varphi^{-1}(U)$ must be isomorphic to the family of curves $\shC'$.
\qed

\begin{corollary}
\mylabel{rational surface}
The surface $G\backslash(C\times C)$ is rational.
\end{corollary}

\proof
The function field of the our quotient surface   is the field of fractions $F$ of the ring
$k[x,y,\tau]/(y^q-y-x^{q+1}-\tau x^q)$.
The residue class of $x$ is clearly a regular element. Thus we have an equation  $\tau = -x+(y^q-y)/x^q$ in the field $F$. Hence
$F=k(x,y)$, which shows that  the surface in question is birational to $\PP^2$.
\qed

\medskip
For $q=2$, Shioda \cite{Shioda 1974} already observed  that the partial quotient by the center 
$Z\backslash(C\times C)$ is rational.
The same thus holds for all further quotients (see for example \cite{Badescu 2001}, Theorem 13.27).

Setting $\tau=1/z$  and clearing denominators for
the equation in   Proposition \ref{twisted form}, we obtain the following description of the generic fiber:

\begin{corollary}
\mylabel{generic equation}
The smooth projective curve $\varphi^{-1}(\eta)$ over the function field  $k(z)$ of the quotient
$\PP^1=G\backslash C$ is given by the homogenization of the equation 
$$
y^q-z^{q^2-1}y=x^{q+1} + z^{q-1}x^q,
$$
where $z\in\O_{\PP^1,\infty}$ is a uniformizer.
\end{corollary}

%===========================================================
\section{Computation of the singular fiber}
\mylabel{Singular fiber}

In this section we shall finish the proof for Proposition \ref{singular fiber}.
In other words, we have to determine the dual graph for the singular fiber 
$\psi^{-1}(\infty)$ on the relative minimal model $\psi:S\ra\PP^1$
of the quotient surface $G\backslash(C\times C)$.
In light of Corollary \ref{generic equation}, the task boils down to compute
the minimal resolution of singularities for the spectrum of residue class ring of $k[x,y,z]$ by the ideal genereated by the equation
\begin{equation}
\label{initial equation}
y^q-z^{q^2-1}y-x^{q+1} - z^{q-1}x^q=0.
\end{equation}
As it turns out, this will produce no $(-1)$-curves in the singular fiber, and thus
yield the desired description of $\psi^{-1}(\infty)$. 
It is worth to examine the simplest special case $q=2$ first. Then $\psi:S\ra\PP^1$ becomes
a jacobian elliptic fibration, whose singular fibers are classified in terms of \emph{Kodaira symbols}.

%-------------------------------------------------------------
\vspace{2em}
\centerline{\includegraphics{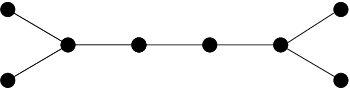}} 
\vspace{1em}
\centerline{Figure \stepcounter{figure}\arabic{figure}: Dual graph for Kodaira symbol ${\rm I}^*_3$}
\vspace{1em}
%-------------------------------------------------------------

\begin{proposition}
\mylabel{singular fiber q=2}
For $q=2$, the Weierstrass equation $y^2-z^3y=x^3+zx^2$ defines a rational double point of type $D_7$,
and the singular fiber of the elliptic fibration $\psi:S\ra\PP^1$ has type ${\rm I}^*_3$.
\end{proposition}

\proof
This follows from the classification of jacobian elliptic fibrations on rational surfaces
obtained by W.\ Lang: Our Weierstrass equation appears as case  in \cite{Lang 1994}, Section 2, Case 13C,
and is listed there with type ${\rm I}^*_3$.
Independently, the structure of the singularity  can be determined with  the 
algorithm of Greuel and Kr\"oning \cite{Greuel; Kroning 1990}.
\qed

\medskip
We now turn to the general case, and depict  the dual graph again in Figure 3,
where  symbols $E_1,\ldots, E_{q-1}, F_0,F_1,F_2,F_3,F_4$ are assigned to certain irreducible components.

%-------------------------------------------------------------
\vspace{2em}
\centerline{\includegraphics{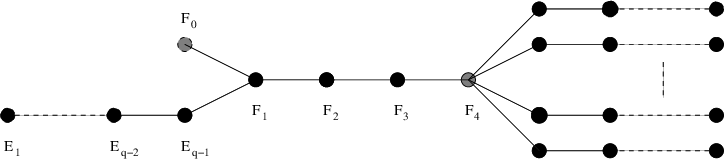}} 
\vspace{1em}
\centerline{Figure \stepcounter{figure}\arabic{figure}: Dual graph of   $\psi^{-1}(\infty)$}
\vspace{1em}
%-------------------------------------------------------------

The resolution of singularities will be carried out in several steps, in which the indicated irreducible components
show up successively. In some sense, our procedure generalizes the part of the \emph{Tate Algorithm} \cite{Tate 1972} dealing with
reduction types ${\rm I}^*_n$, $n\geq 0$. The  main idea already present in the Tate Algorithm
is to blow-up Weil divisors, rather than   points.

\medskip
{\bf Step 0.} The fiber for the initial equation (\ref{initial equation}) is given by $z=0$,
hence equals the spectrum of the ring $k[x,y]/(y^q-x^{q+1})$, which is clearly a rational curve
with a single singularity at the origin. It corresponds to the irreducible component $E_1$ in Figure 3.

\medskip 
{\bf Step 1.} Let us pass to the variables $z,x/z^{q-1},y/z^{q-1}$, and examine the strict transform
\begin{equation}
\label{weighted blowup}
\left(\frac{y}{z^{q-1}}\right)^q - z^{2q-2}\left(\frac{y}{z^{q-1}}\right) - \left(\frac{x}{z^{q-1}}\right)^{q+1}z^{q-1} - z^{q-1}\left(\frac{x}{z^{q-1}}\right)^q = 0
\end{equation}
of our initial equation (\ref{initial equation}), which is obtained by making   substitutions
$x=x/z^{q-1}\cdot z^{q-1}$ and  $y=y/z^{q-1}\cdot z^{q-1}$,
and afterwards dividing by the highest possible $z$-power.
The schematic fiber on this chart is given by $z=0$, and   thus a copy of the affine line with multiplicity $q$.
This will be the irreducible component $F_1$ in Figure 3. Note that for $q\neq 2$, the   Equation (\ref{weighted blowup})
defines a divisor that is singular along the fiber.

\medskip 
{\bf Step 2.} We next blow-up the reduced fiber for the previous Equation (\ref{weighted blowup}). 
In other words, the center is given by the ideal $(z,y/z^{q-1})$.
This gives us two new charts. We start by describing the $y/z^{q-1}$-chart, which features as new variables 
$z^q/y,x/z^{q-1},y/z^{q-1}$,
allowing for the substitution $z=z^q/y\cdot y/z^{q-1}$. The resulting strict transform is
\begin{equation}
\label{}
\frac{y}{z^{q-1}} - \left(\frac{z^q}{y}\right)^{2q-2} \left(\frac{y}{z^{q-1}}\right)^q 
- \left(\frac{x}{z^{q-1}}\right)^{q+1} \left(\frac{z^q}{y}\right)^{q-1} - \left(\frac{z^q}{y}\right)^{q-1}\left(\frac{x}{z^{q-1}}\right)^q=0.
\end{equation}
Computing partial derivatives one sees that this chart is smooth.
The fiber is given by $z^q/y\cdot y/z^{q-1}=0$, thus becomes the spectrum of residue class ring of $k[z^q/y,x/z^{q-1}]$ 
by the equation 
$$
\left(\frac{z^q}{y}\right)^q \left(\frac{x}{z^{q-1}}\right)^q \left(\frac{x}{z^{q-1}}-1\right)=0.
$$
The three factors give   the three irreducible components $qF_1+qF_2+F_0$ in Figure 3, respectively.
To see that there are no further intersections, one has to look at the other chart, the $z$-chart.
Here the variables are $z,x/z^{q-1},y/z^q$. The substitution is $y/z^{q-1}=y/z^q\cdot z$, and the
resulting strict transform is
\begin{equation}
\label{z-chart}
\left(\frac{y}{z^q}\right)^q z - z^q \frac{y}{z^q} - \left(\frac{x}{z^{q-1}}\right)^{q-1} - \left(\frac{x}{z^{q-1}}\right)^q = 0.
\end{equation}
The fiber is given by $z=0$, whence isomorphic to $k[x/z^{q-1},y/z^q]$ modulo the equation 
$$
\left(\frac{x}{z^{q-1}}\right)^q\left(\frac{x}{z^{q-1}}-1\right).
$$
This shows that there are no further intersection numbers. Computing partial derivatives,
one sees that the   Equation (\ref{z-chart}) has just one singularity, which is located at the origin.

\medskip
{\bf Step 3.}  Now we make a blowing-up of (\ref{z-chart}) along the reduced fiber, which is given by the ideal
$(z,x/z^{q-1})$.
Let us treat the $z$-chart first. Here we have variables $z,x/z^q,y/z^q$, resulting in the  strict transform
\begin{equation}
\label{second z-chart}
\left(\frac{y}{z^q}\right)^q - z^{q-1}\frac{y}{z^q} - \left(\frac{x}{z^q}\right)^{q+1}z^q - z^{q-1} \left(\frac{x}{z^q}\right)^q = 0.
\end{equation}
The fiber on this chart is given by $z=0$, which is a copy of the affine line with multiplicity $q$.
It corresponds to the irreducible component $qF_3$ in Figure 3. Computing partial derivatives, one sees that
the singular locus is contained in the fiber.
A computation of the $x/z^{q-1}$-chart left to the reader reveals that $F_2\cdot F_3=1$, and that there are no
additional singularities.

\medskip
{\bf Step 4.} We  blow-up the previous equation (\ref{second z-chart}) along the  reduced fiber for the previous equation, 
which is given by the ideal $(z,y/z^q)$.
The $z$-chart has variables $z,x/z^q,y/z^{q+1}$, and the ensuing strict transform is
\begin{equation}
\label{}
\left(\frac{y}{z^{q+1}}\right)^q z - z\frac{y}{z^{q+1}} - \left(\frac{x}{z^q}\right)^{q+1} z -  \left(\frac{x}{z^q}\right)^q = 0.
\end{equation}
The fiber on this chart is given by $z=0$, which is a copy of the affine line with multiplicity $q$.
It corresponds to the irreducible component $F_4$ in Figure 3. Computing partial derivatives, one sees that the singular locus is
given by $z=x/z^q=0$ and $(y/z^{q+1})^q-y/z^{q+1}=0$. In other words, there are exactly $q$ singular
points, in a canonical way  parameterized by the finite field $\FF_q$. Looking at the formal completions
of the local rings at these singularities, one easily sees that the singular points are rational double points
of type $A_{q-1}$, that is, formally isomorphic to $k[[u,v,w]]/(u^q-vw)$, and that the resulting irreducible components
intersect $F_4$ as in Figure 3. This accounts for the $q$ strings of length $q-1$ on the right hand side of Figure 3.
Computing the $y/z^q$-chart, one checks that $F_3\cdot F_4=1$. 

\medskip
{\bf Step 5.} It remains to understand how the string  of $q-1$ curves on the left hand side of Figure 3
arises. For this return to the initial equation (\ref{initial equation}), and blow-up 
the \emph{nonreduced 0-dimensional center} given by $(z^{q-1},x,y)$. We have already computed the $z$-chart in Step 1.
The task now is to understand the $x$-chart. Here we have four variables $z^{q-1}/x,z,x,y/x$, and the strict transform is described
by   two equations
\begin{equation}
z^{q-1} = z^{q-1}/x\cdot x\quadand 
\left(\frac{y}{x}\right)^q - \left(\frac{z^{q-1}}{x}\right)^2 \frac{y}{x} x^2 - x - \frac{z^{q-1}}{x}x = 0.
\end{equation}
The fiber $z=0$ becomes thus the spectrum of the residue class ring of  $k[z^{q-1}/x,x]$ by the equation
$z^{q-1}/x\cdot(\frac{y}{x})^q=0$. Its two irreducible components
correspond to $E_1$ and $qF_1$. 
Now we complete the local ring at the origin. The second equation, together with the Implicit Function Theorem, 
tells us that $x=u(y/x)^q$, where $u$ is in invertible formal power series. By the first equation,
our complete local ring is isomorphic
to 
$$
k[[u,v,w]]/(u^{q-1}-vw^q).
$$
By the classical theory of quotient singularities, the normalization of this ring contains a unique singularity,
which is a Hirzebruch--Jung singularity, compare \cite{Barth; Peters; Van de Ven 1984}, Chapter III, Section 5.
The self-intersection numbers of the irreducible components on the
minimal resolution of singularities are given by the continued fraction development:
$$
\frac{q-1}{q-2} = 
\underbrace{
2- \cfrac {1}{\ddots-
    \cfrac{1}{2}}}_{\text{$q-2$ entries}}
$$
In other words, we have a string of $q-2$ rational curves with self-intersection $(-2)$, so our singularity
is actually a rational double point of type $A_{q-2}$. This gives the irreducible components $E_2,\ldots,E_{q-1}$
in Figure 3. One easily computes the  multiplicities in the fibers are $\sum_{i=2}^{q-1}iE_i$.
A straightforward computation of the $y$-chart left to the reader reveals that there are no further singularities appear.

\medskip
{\bf Step 6.}
By now we have verified that the dual graph of the fiber $\psi^{-1}(\infty)$ is as in Figure 1, and
that the multiplicities of the irreducible components are as indicated. Using the fact that the divisor $\psi^{-1}(\infty)$
is numerically trivial on itself, one easily computes the self-intersection numbers by induction, starting with the irreducible
components having only one neighbor.
This concludes the proof for Proposition \ref{singular fiber}, and thus also for Theorem \ref{two nodes}.

%===========================================================
\section{Higher ramification groups and Swan conductor}
\mylabel{Higher Ramification}

We now take a closer look at the action of the special $p$-group $G$ on the Hermitian curve
$C:\; y^q-q=x^{q+1}$ at infinity.
The closure of the affine part of the scheme $C$ inside $\PP^2$ is given by the homogeneous 
equation
$$
Y^qZ-YZ^q =X^{q+1},\qquad y=Y/Z,\quad x=X/Z.
$$
Clearly, only the point $(0:1:0)\in\PP^2$ lies  on the closure at infinity. Dehomogenizing with
respect to $Y$, we obtain the equation
$$
Z/Y - (Z/Y)^q = (X/Y)^{q+1}.
$$
Completing and using the implicit function theorem, one sees that 
$$
k[[X/Y]] = k[[X/Y,Z/Y]]/(Z/Y - (Z/Y)^q - (X/Y)^{q+1}),
$$
and furthermore
$$
Z/Y = (X/Y)^{q+1} + \ldots.
$$
Note that here and throughout,  dots in power series expansions are shorthand notation for ``higher order terms''.

Summing up, we have an embedding $C\subset\PP^2$ of the regular proper curve $C$, and
the complete local ring of the curve at the point $c=(0:1:0)$  is
$\O_{C,c}^\wedge = k[[X/Y]]$.
The action of $G$ on $C$ leaves $c\in C$ fixed. We now compute its \emph{higher ramification
groups} 
$$
G_i=\left\{\sigma\in G\mid \text{$\sigma=\id$ on $\O_{C,c}^\wedge/\maxid^{i+1}$}\right\}.
$$
These form a descending filtration $G_0\supset G_1\supset G_2\supset\ldots$ 
of normal subgroups inside $G$.

\begin{proposition}
\mylabel{higher ramification}
The higher ramification groups are 
$$
G_0=G_1=G,\quad  G_2=\ldots=G_{q+1}=Z,\quad G_{q+2}=0.
$$
\end{proposition}

\proof
We first determine  
the position in the descending filtration for the nontrivial central group elements  $\sigma=(t,0)$, $t\neq 0$. 
The action  on the uniformizer in $\O_{C,c}^\wedge$ is
$$
X/Y=x/y \longmapsto x/(y+t) = x/y \cdot (1+t/y)^{-1}.
$$
We have $1/y=Z/Y=(X/Y)^{q+1}+\ldots$, whence 
$$
(1+t/y)^{-1}=\sum_{n\geq 0} (-t/y)^n = 1 - t(X/Y)^{q+1} + \ldots\;.
$$
The upshot is that $\sigma\in G_{q+1}$ but $\sigma\not\in G_{q+2}$.
Now we consider noncentral group elements $\sigma=(t,r)$, $r\neq 0$. Here the
action on the uniformizer is 
$$
X/Y=x/y\longmapsto (x+r)/(y-r^q x + t) = x/y \cdot (1+r/x) \cdot (1-r^q x/y + t/y)^{-1}.
$$
The factors on the right hand side come from
$$
1-r^q x/y +t/y = 1- r^q(X/Y) + t(X/Y)^{q+1} + \ldots
$$
and
$$
1+r/x=1+r\cdot Z/Y\cdot (X/Y)^{-1}  = 1 + r(X/Y)^q + \ldots\; .
$$
Arguing as in the preceding paragraph with the geometric series, we infer that
$\sigma \in G_1$ but $\sigma\not\in G_2$.
\qed

\medskip
Now let $l$ be a prime   different from $p$, and let $M$ be a finite dimensional $\FF_l$-vector space
endowed with a $G$-representation. Then one way to define the \emph{Swan conductor} is
$$
\delta=\delta(G,M) = \sum_{i\geq 1} \frac{1}{[G:G_i]} \dim(M/M^{G_i}).
$$
As explained in \cite{Serre 1978}, Chapter 19, this rational number is actually an integer. For later
use we record:

\begin{corollary}
\mylabel{swan conductor}
In our situation, the Swan conductor is
$$
\delta=\dim(M/M^G) + q^{-1}\dim(M/M^Z).
$$
\end{corollary}

%===========================================================
\section{Representations of certain special $p$-groups}
\mylabel{Representations}

In the next section we shall determine the $l$-adic cohomology group $H^1(C,\QQ_l)$ as $G$-representation, where
$C:\;y^q-y=x^{q+1}$ is a Hermitian curve, and $G$ is the special $p$-group acting on $C$.
To this end, we now make a digression and examine the representation theory of $G$ over various fields. 
For general facts about representations of finite groups, we refer to the monographs of
Serre \cite{Serre 1978},  Curtis and Reiner \cite{Curtis; Reiner 1990} and Isaacs \cite{Isaacs 1976}.
The representation theory  of extraspecial $p$-groups is well-known (for example \cite{Doerk; Hawkes 1992}, Section B.9
or \cite{James; Liebeck 2001}, Section 26).

It   seems natural to deal with a somewhat more general class of groups.
Throughout this section, we fix a prime power $q=p^n$, and 
assume that $G$ is a finite special $p$-group. To keep the exposition in bounds,  we demand that $\exp(G)=p$ for odd $p$.
In contrast, we stipulate that   each noncentral element has order four in case $p=2$, as suggested by   Proposition \ref{exponent}.

Moreover, we suppose that the group is
endowed with the following \emph{additional structure}:
The center $Z$ and the abelianization $G^\ab=G/Z$, which are elementary abelian groups and thus   $\FF_p$-vector spaces, 
are endowed with an $\FF_q$-vector space structure so that the commutator pairing
$$
G^\ab\times G^\ab\lra Z,\quad (\bar{a},\bar{b})\longmapsto [a,b]=aba^{-1}b^{-1}
$$
are $\FF_q$-bilinear and $\dim_{\FF_q}(Z)=1$.
In turn, the abelianization $G^\ab$ may be regarded as a symplectic $\FF_q$-vector space.
Obviously, the order of the group is of the form 
$$
\ord(G)=q^{1+2n},\qquad n=\dim_{\FF_q}(G^\ab).
$$
Clearly, the automorphism group $G\subset\Aut(C)$ considered in Section \ref{Hermitian curves} is   endowed with such an
additional structure in an obvious way. To start with, we record:

\begin{proposition}
\mylabel{conjugacy classes}
The conjugacy classes of the special $p$-group $G$ are the subsets
$$
\left\{z\right\}\subset G \quadand f^{-1}(\bar{b})\subset G,
$$
with $z\in Z$ and $\bar{b}\in G^\ab\smallsetminus 0$, where $f:G\ra G^\ab$ is the canonical projection.
\end{proposition}

\proof
We only have to verify that the fibers $f^{-1}(\bar{b})$, $\bar{b}\neq 0$ are whole conjugacy classes. Since $b\notin Z$, there
is a nontrivial commutator $[a,b]$. Since the commutator pairing is $\FF_q$-linear in $\bar{a}$, the commutators $[a,b]$, $a\in G$
fill out the whole center $Z\subset G$. In turn, 
the conjugate elements $aba^{-1}=[a,b]\cdot b$ fill out the whole coset $f^{-1}(\bar{b})\subset G$.
\qed

\medskip
Given  an arbitrary field $K$, we write  $R_K(G)$  for the Grothendieck group of finite-dimensional
$G$-representations over $K$, modulo the relations $[V]=[V']+[V'']$ coming from short exact sequences
$0\ra V'\ra V\ra V''\ra 0$.  Let
$$
\Irr_K(G)=\left\{[V_1],\ldots,[V_r]\right\}
$$
be the set of isomorphism classes of irreducible   $G$-representations over $K$ or, what amounts to the same, simple $KG$-modules.
Then $\Irr_K(G)$ is a basis for the abelian group $R_K(G)$, and we say that
the  $V_1,\ldots,V_r$ from a \emph{basic set of irreducible representations}. 
We now shall determine basic sets of irreducible $G$-representations for various fields $K$.
Our further analysis   hinges on the following observation:

\begin{proposition}
\mylabel{irreducible representations}
Let $K$ be a field of characteristic zero. Then  the set $\Irr_K(G)$ has cardinality 
$$
\Card(\Irr_K(G))=
\begin{cases}
1+(q^{2n}+q-2)/(p-1)  & \text{if $\mu_p(K)=\left\{1\right\}$;}\\
q^{2n}+q-1            & \text{else.}
\end{cases}
$$
\end{proposition}

\proof
First note that over algebraically closed fields of characteristic zero, the number of irreducible representations
equals the number of conjugacy classes of elements; then our statement follows from Proposition \ref{conjugacy classes}.

Now let $\xi\in\mu_{p^2}(\bar{K})$ be a primitive $p^2$-th root of unity. The field extension $E=K(\xi)$ is 
abelian, with Galois group $\Gamma\subset(\ZZ/p^2\ZZ)^\times$.
Since each element $a\in G$ has $a^{p^2}=e$, we have a natural action of the group $\Gamma$ on the set $G$,
where $\sigma\in\Gamma$ acts via $a\mapsto a^\sigma$. According to \cite{Serre 1978}, Section 12.4, the
cardinality of $\Irr_K(G)$ coincides with the number of $\Gamma$-orbits of conjugacy classes in $G$.
Since $Z$ and $G^\ab$ are elementary abelian, 
the  $p$-torsion subgroup $(1+p\ZZ)/p^2\ZZ\subset\Gamma$ acts trivially on the
the set of conjugacy classes, and the formula in the second case of the statement follows.

Now suppose that $\mu_p(K)=\left\{1\right\}$. Then $\Gamma$ acts via its quotient $(\ZZ/p\ZZ)^\times$, and the latter group 
permutes transitively the nonzero central elements $z\in Z$ generating the same subgroup, and permutes transitively the 
fibers $f^{-1}(\bar{b})$ where $\bar{b}\in G^\ab$, $\bar{b}\neq 0$   generates the same subgroup in $G^\ab$. 
Clearly,   center and abelianization contain $(q-1)/(p-1)$ and $(q^{2n}-1)/(p-1)$ nonzero cyclic subgroups, respectively,
and  the first case of the statement follows.
\qed

\medskip
Our next task is to construct and describe the irreducible $G$-representations over $K=\QQ$.
Of course there is the trivial representation, which we simply denote by $\QQ$.
Next, let $\varphi:G\ra\CC^\times$ be a nonzero homomorphism; in other words, a   nonzero homomorphism
$\varphi^\ab:G^\ab\ra\CC^\times$.   It factors over the subgroups $\mu_p(\CC)\subset\QQ(e^{2\pi i/p})^\times\subset\CC^\times$,
because $G^\ab$ is elementary abelian,
so we   obtain a $G$-representations $(\QQ(e^{2\pi i/p}),\varphi)$ over $\QQ$ of degree $p-1$.
Up to isomorphism, this representation depends only on the $\FF_p$-hyperplane $\ker(\varphi)\subset G^\ab$,
as the action of the Galois group for the Galois extension $\QQ\subset\QQ(e^{2\pi i/p})$ shows.
Using Grothendieck's convention, we  regard such hyperplanes as rational points $y\in\PP(G^\ab)\simeq\PP^{mn-1}_{\FF_p}$,
and set
$$
W_y = (\QQ(e^{2\pi i/p}),\varphi).
$$
By the very construction, its degree and isomorphism algebra is given by
\begin{equation}
\label{info W}
\dim(W_y)=p-1\quadand \End(W_y)=\QQ(e^{2\pi i/p}).
\end{equation}
In particular, the representation is irreducible, and its complexification splits as $W_y\otimes_\QQ\CC =\bigoplus(\CC,\varphi)$, where the sum
extends over all nonzero homomorphisms $\varphi:G\ra\CC^\times$ satisfying $y=\ker(\varphi^\ab)$.

To construct the remaining irreducible representations, one makes a   construction involving induction from the center $Z\subset G$.
Choose once and for all a
Lagrangian $\FF_q$-subvector space in $G^\ab$, that is, a totally isotropic vector subspace of maximal dimension,
and let $H\subset G$ be its preimage. This is an abelian normal subgroup of
index $[G:H]=q^n$ with $Z\subset H$.
By our standing assumption $Z\subset H$ is a direct summand if and only if $p$ is odd.
At this point it is   most convenient to distinguish that cases that $p$ is even or odd.

Suppose $p$ is odd.
Given a nonzero homomorphism $\chi:Z\ra\CC^\times$, we extend it 
in an arbitrary way to a homomorphism $\tilde{\chi}:H\ra\CC^\times$. Since $H$ is elementary abelian, it factors
over the subgroups $\mu_p(\CC)\subset\QQ(e^{2\pi i/p})\subset\CC^\times$, thus turning $\QQ(e^{2\pi i/p})$ into an irreducible
$H$-representation of degree $p-1$. As above, the Galois action shows that the underlying $Z$-representation representation depends only 
$\ker(\chi)\subset Z$, which we regard as a rational point $x\in\PP(Z)\simeq\PP^{m-1}_{\FF_p}$. Now consider the induced $G$-representation
$$
V_x=\ind_H^G(\QQ(e^{2\pi i/p}),\tilde{\chi})=\QQ(e^{2\pi i/p})\otimes_{\QQ H}\QQ G 
$$
over $\QQ$. A priori, this seems to depend on the chosen extension $\tilde{\chi}$. However, it will follow from
Theorem \ref{basic p>2} below that its isomorphism
class depends only on $x$.
																	     
\begin{lemma}
\mylabel{irreducible p>2}
Suppose $p\neq 2$. The $G$-representations $V_x$, $x\in\PP(Z)$  are irreducible. Their degrees and endomorphism algebras are
$$
\dim(V_x)=q^n(p-1)\quadand \End(V_x)=\QQ(e^{2\pi i/p}).
$$ 
\end{lemma}

\proof
The degree is obvious, and   irreducibility follows from the statement about the endomorphism ring.
By definition, we have $V_x=\ind^G_H(M)$ for the irreducible $H$-represen\-tation $M=(\QQ(e^{2\pi i/p}),\tilde{\chi})$.
Frobenius Reciprocity gives
$$
\End(V_x)=\Hom_H(M, \res^G_H(\ind_H^G(M))),
$$
and $\res^G_H(\ind_H^G(M))$ is the direct sum of the $H$-representations ${}^bM$ obtained by transport of structure, 
where $bH\subset G$ runs over all cosets.
The $H$-action on ${}^bM$ is defined as $bhb^{-1} \cdot m=hm$.
Therefore, we merely
have to check that $\chi\circ\gamma_b\neq\chi$ for all $b\in G\smallsetminus H$.  Since Lagrangian are maximally totally isotropic, there must be some $a\in H$ that does
not commute with $b$. Now choose some some $z\in Z$ not contained in
the kernel of $\chi:Z\ra\CC^\times$.  Multiplying by a suitable $\lambda\in\FF_q$ if necessary, we may assume that $[a,b]=z$.
Then $\chi(bab^{-1})\neq\chi(a)$.
\qed

\begin{theorem}
\mylabel{basic p>2}
Suppose $p\neq 2$. Then the representations
$$
V_x,\, x\in\PP(Z),\quadand W_y,\, y\in\PP(G^\ab),\quadand \QQ
$$
form a basic set of irreducible representations over $K=\QQ$. Moreover, they remain a basic set of irreducible representations
after scalar extensions $\QQ\subset K$ if and only if  $\mu_p(K)=\left\{1\right\}$.
\end{theorem}

\proof
We saw above that the given representations are irreducible. It is obvious that the $(q^{2m}-1)/(p-1)$ representations
$W_y$ are pairwise nonisomorphic,
and that no $V_x$ is isomorphic to some $W_y$. To see that the $(q-1)/(p-1)$ representations $V_x$ are pairwise nonisomorphic,    look
at their characters on $Z$. For $K=\QQ$, the statement now follows from Proposition \ref{irreducible representations}.

Now let $\QQ\subset K$ be a field extension with $\mu_p(K)=\left\{1\right\}$. Then the endomorphism algebra of
$V_x\otimes_\QQ K$ is isomorphic to $\QQ(e^{2\pi i/p})\otimes_\QQ K$, which remains a field, so the representation stays
irreducible. The same argument applies for the representations $W_y$ and $\QQ$. A second application of Proposition \ref{irreducible representations} concludes the proof.
\qed

\medskip
Now we have to deal with the case $p=2$. Let $\chi:Z\ra\CC^\times$ be a nonzero homomorphism, and extend it
in an arbitrary way to $\tilde{\chi}:H\ra \CC^\times$. Then $\chi$ factors over $\mu_2(\CC)=\left\{\pm 1\right\}$,
yet this is not the case for $\tilde{\chi}$, because $Z\subset H$ is not a direct summand by our standing assumption.
However, it factors over $\mu_4(\CC)\subset\QQ(i)^\times\subset\CC^\times$. We now define a $G$-representation
$$
2V_x = \ind_H^G(\QQ(i),\tilde{\chi}) = \QQ(i)\otimes_{\QQ H}\QQ G.
$$
Here $x\in\PP(Z)$ is the rational point corresponding to $\chi:Z\ra\CC^\times$ via $x=\ker(\chi)$.  
Beware that  $2V_x$ is just a formal symbol; there is no representation $V_x$ over $\QQ$; the latter
will acquire meaning only after certain scalar extensions.

Now recall that   $(a,b)_K$ with $a,b\in K^\times$ denotes the \emph{quaternion algebra} over a field $K$ generated by
symbols $i,j$ modulo the relations $i^2=a$, $j^2=b$ and $ij=-ji$. This a central simple algebra of degree two. Up to isomorphism,
it depends only on the classes of $a,b$ modulo $\QQ^{\times 2}$. For example,
$(-1,-1)_\RR=\HH$ are the classical Hamilton quaternions.

\begin{lemma}
\mylabel{irreducible p=2}
Suppose $p=2$. The $G$-representations $2V_x$, $x\in\PP(Z)$ are irreducible. Their degree and endomorphism algebra is given
by
$$
\dim(2V_x)=2q^n\quadand \End(2V_x)=(-1,-1)_\QQ.
$$ 
\end{lemma}

\proof
It is easy to see that $\QQ(i)$ is irreducible as $Z$-representations, and the irreducibility of the induced
$G$-representation follows with Mackey's criterion as in the proof of Lemma \ref{irreducible p>2}.
The dimension formula is obvious. 

The endomorphism algebra $D=\End(2V_x)$ is a division ring.
According to Proposition \ref{irreducible representations}, the cardinality of $\Irr_K(G)$ does not depend on the field $\QQ\subset K$.
We infer that $\QQ$ must be the center of the division ring $D$. 
According to Roquette's result (compare \cite{Isaacs 1976}, Corollary 10.14), the degree $d=\deg(D)$, which is defined as the square root of $\dim_\QQ(D)$ and
in this context called the \emph{Schur index}, is either $d=1$ or $d=2$. We have an obvious inclusion $\QQ(i)\subset D$.
Whence $d=2$, and our division ring is of the form $D=(a,b)_\QQ$, and we may assume that $a,b\neq 0$ are square-free integers.
Clearly $K=\RR$ is not a splitting field for the $H$-representation $\QQ(i)$, and therefore not for $2V_x$.
It follows $a,b<0$. Similarly, $K=\QQ(\sqrt{n})$ is not a splitting field, for all square free integers $n>0$.
So the only remaining possibility is $a,b=-1$.
\qed

\begin{theorem}
\mylabel{basic p=2}
Suppose $p= 2$. Then the representations
$$
2V_{x},\, x\in\PP(Z)\quadand W_{y},\, y\in\PP(G^\ab)\quadand \QQ
$$
form a basic set of irreducible representations over $K=\QQ$. Moreover, they
stay a basic set of irreducible representations after scalar extension $\QQ\subset K$
if and only if the equation $T_0^2+T_1^2+T_2^2=0$ has no nontrivial solution in $K$.
\end{theorem}

\proof
We saw above that the given representations are irreducible. It is obvious that the $(q^{2m}-1)/(p-1)$ representations
$W_y$ are pairwise nonisomorphic,
and that no $V_x$ is isomorphic to some $W_y$. To see that the $(q-1)/(p-1)$ representations $V_x$ are pairwise nonisomorphic, simply  look
at their characters on $Z$. For $K=\QQ$, the statement now follows from Proposition \ref{irreducible representations}.
For $\QQ\subset K$ arbitrary, $2V_x$ stays irreducible over $K$ is and only if this field does not split
$\End(2V_x) =(-1,-1)_\QQ$, whence the result.
\qed

\medskip
If $\QQ\subset K$ is any field extension, we set
$$
V_{x,K}= V_x\otimes_\QQ K \quadand W_{y,K}=W_y\otimes_\QQ K
$$
for the $G$-representations obtained by scalar extension. 
Note that for $p=2$, the former only makes sense if $T_0^2+T_1^2+T_2^2=0$ has a nontrivial solution in $K$;
in this case, $V_{x,K}$ is defined as an irreducible summand in $2V_x\otimes_\QQ K$.

Over the field $K=\QQ_l$ of $l$-adic numbers, $l\neq p$, we have the following almost uniform description of irreducible
representation:

\begin{proposition}
\mylabel{basic set}
Suppose either $p=2$, or  that $p$ does not divide $l-1$. Then the 
$V_{x,\QQ_l},\, W_{y,\QQ_l},\, \QQ_l$ form a basic set of irreducible $G$-representations over $K=\QQ_l$.
\end{proposition}

\proof
Suppose first that $p$ is odd. By assumption, $\mu_p(\FF_l)$ is trivial, thus the same holds for $\mu_p(\QQ_l)$,
and the statement follows from Proposition \ref{basic p>2}.
Now suppose $p=2$. The quadric $T_0^2+T_1^2+T_2^2$ has a nontrivial solution over any finite field, thus also over $\QQ_l$ by
Hensel's Lemma, so $V_{x,\QQ_l}$ makes sense. By (\ref{info W}), the $W_{x,\QQ_l}$ stay irreducible, and
the statement easily follows from Proposition \ref{basic p=2}.
\qed

%===========================================================
\section{Cohomology as representations}
\mylabel{Cohomology as}

We return to our  Hermitian curve
$C:\; y^q - y = x^{q+1}$, on which the special $p$-group $G$ of order $\ord(G)=q^3$ acts,
as discussed in Section \ref{Hermitian curves}.  
We seek to compute the $l$-adic cohomology group $H^1(C,\QQ_l)$, viewed as a $G$-representation of degree $q(q-1)$ over the field $K=\QQ_l$.
In the previous section, we have introduced the   $G$-representation
$V_{x,\QQ_l}=V_x\otimes_\QQ\QQ_l$. We now set
$$
V_{\QQ_l} =\bigoplus_{x} V_{x,\QQ_l}
$$
where the sum runs over all points $x\in\PP(Z)\simeq\PP^{m-1}_{\FF_p}$. Recall that $Z\subset G$ is the center.

\begin{theorem}
\mylabel{cohomology representation}
For each prime $l\neq p$, the $G$-representations  $H^1(C,\QQ_l)$ and $V_{\QQ_l}$ are isomorphic.
\end{theorem}

\proof
Choose a nonzero element $z\in Z$, and consider its  \emph{Lefschetz number}
$$
\Lambda(z) = \sum_{i=0}^2 (-1)^i\Tr(z, H^i(C,\QQ_l))\in\QQ_l.
$$
By the Lefschetz Trace  Formula (for example \cite{Milne 1980}, Theorem 12.3), the number $\Lambda(z)$ equals the length of the fixed scheme on $C$.
We saw in Proposition \ref{higher ramification} that $z\in G_{q+1}$ but $z\notin G_{q+2}$, hence the fixed scheme
has length $q+2$. We conclude
$$
\Tr(z,H^1(C,\QQ_l)) = -\Lambda(z) + 2 = -q.
$$
Before we continue, note that it suffices to argue for a single prime $l\neq p$, because
the irreducible components in all the  $V_{x,\CC}$ are pairwise nonisomorphic.
So we may assume that $p$ does not divide $l$ in case $p\neq 2$.
Now we are in a situation that the $V_{x,\QQ_l},\, W_{y,\QQ_l},\, \QQ_l$ form a basic set of irreducible representations
over $K=\QQ_l$. Write the isomorphism class of $H^1(C,\QQ_l)$ as
$$
\sum_{x\in\PP(Z)} m_x[V_{x,\QQ_l}] + \sum_{y\in\PP(G^\ab)} n_y [W_{y,\QQ_l}] + n[\QQ_l],
$$
for certain integral coefficients $m_x,n_y,n\geq 0$. Obviously, $\Tr(z,W_y)=p-1$ and $\Tr(z,\QQ)=1$
are positive. As the Lefschetz number $\Lambda(z)$ is negative, at least one $m_x$ must be nonzero.
Since
$$
\dim H^1(C,\QQ_l)=q(q-1) = \frac{q-1}{p-1}q(p-1) = \dim(V_{\QQ_l}),
$$
it suffices to verify that the coefficients $m_x$ do not depend on the index $x\in \PP(Z)$.
To achieve this, we use   further symmetries of the Hermitian curve $C:\;y^q-y=x^{q+1}$.
Namely, the group $\mu_{q^2-1}(k)=\FF_{q^2}^\times$ acts on $C$ via
$$
x\longmapsto \zeta x,\quad y\longmapsto \zeta^{q+1}y,\qquad \text{for $\zeta\in\FF_{q^2}^\times$,}
$$
as remarked in \cite{Stichtenoth 1973b}, proof of Satz 5.
One easily computes that the subgroup $\FF_{q^2}^\times\subset\Aut(C)$ normalizes the subgroup
$G\subset\Aut(C)$, and the induced conjugacy action
$$
c: \FF_{q^2}^\times\lra \Aut(G),\quad \zeta\longmapsto \left( (t,r)\mapsto \zeta\circ(t,r)\circ\zeta^{-1}\right)
$$
is given by $(t,r)\mapsto(\zeta^{q+1}t,\zeta r)$. We   infer that  the   linear bijection $\zeta:H^1(C,\QQ_l)\ra H^1(C,\QQ_l)$
is an isomorphism between the original $G$-representation and the new $G$-representation obtained by transport
of structure
$$
G\stackrel{c_\zeta}{\lra} G\lra \GL (H^1(C,\QQ_l)).
$$
A straightforward computation shows that
the $G$-representation ${}^\zeta V_x$ obtained by transport of structure is precisely $V_{\zeta^{q+1}x}$.
But  the $\zeta^{q+1}$ range  over all elements in $\mu_{q-1}(k)=\FF_q^\times$; thus the group $\FF_{q^2}^\times$ acts transitively
on the    nonzero elements in $Z$, and in particular transitively on the $\FF_p$-hyperplanes
$x\subset Z$. From this we deduce that the multiplicities $m_x\geq 0$ do not depend on the index $x\in\PP(Z)$.
\qed

\medskip
For later use, we record the following consequence:

\begin{corollary}
\mylabel{invariant cohomology}
For each prime $l\neq p$, the $G$-invariant part of $H^1(C,\QQ_l)$ is trivial, whereas the $G$-invariant part of
the tensor product $H^1(C,\QQ_l)\otimes_{\QQ_l} H^1(C,\QQ_l)$ has dimension $q-1$.
\end{corollary}

\proof
As in the previous proof, it suffices to verify this for a single prime $l\neq p$, and we thus may assume that
the $V_{x,\QQ_l}$ are irreducible. The first statement immediately follows from the Theorem.
As to the second part, it suffices to check that the $G$-invariant part of $V_{x,\QQ_l}\otimes V_{x,\QQ_l}$ has
dimension $p-1$, because $H^1(C,\QQ_l)$ contains $(q-1)/(p-1)$ such irreducible summands.
Looking at traces, one easily sees that the $V_{x,\QQ_l}$ are selfdual, hence
we have
$$
V_{x,\QQ_l}\otimes V_{x,\QQ_l}\simeq\End(V_{x,\QQ_l}).
$$
The statement follows from Lemmas \ref{irreducible p>2} and \ref{irreducible p=2}.
\qed

\begin{corollary}
\mylabel{cohomology irreducible}
Let $l\neq p$ be a prime. Suppose either $p=2$, or that $p$ does not divide $l-1$.
Then $H^1(C,\QQ_l)$ is irreducible as representation of $\Aut(C)$.
\end{corollary}

\proof
According to the Theorem, we have a decomposition
$$
H^1(C,\QQ_l)=\bigoplus_{x\in\PP(Z)} V_{x,\QQ_l}.
$$
The assumption ensures that  the summands $V_{x,\QQ_l}$ are irreducible by Proposition \ref{basic set}.
In the proof of the Theorem, we saw that the summands are permuted transitively
by the subgroup $\FF_{q^2}^\times\subset\Aut(C)$. It follows that $H^1(C,\QQ_l)$
is irreducible as representation of $G\rtimes\FF_{q^2}^\times\subset\Aut(C)$.
\qed

\medskip
Applying Corollary \ref{swan conductor}, one easily obtains:

\begin{corollary}
\mylabel{wild conductor}
Let $l\neq p$ be a prime. The Swan conductor for the $G$-representation $H^1(C,\FF_l)$
is given by $\delta(G,H^1(C,\FF_l))=q^2-1$.
\end{corollary}

%===========================================================
\section{Curves with irreducible cohomology}
\mylabel{Curves with}

In this section we shall examine curves  whose first $l$-adic cohomology   is 
irreducible as representation of the full automorphism group 
for some prime $l\neq p$. The precise setting is as follows.
Fix a prime power $q=p^m$, and consider smooth proper geometrical connected curves $C_0$
over the finite field $\FF_q$, and let  $C=C_0\otimes_{\FF_q} k$ be
the induced curve over the algebraic closure $k=\bar{\FF}_q$.

The   power $\Fr_{C_0}^m$ of the absolute Frobenius morphism $\Fr_{C_0}:C_0\ra C_0$  is a $\FF_q$-morphism,
and induces by base change a bijective morphism
$$
\Phi=\Fr_{C_0}^m\otimes_{\FF_q} k:C\lra C.
$$ 
We now consider the induced action 
$$
\Phi: H^1(C,\QQ_l)\lra H^1(C,\QQ_l),\quad e\geq 1
$$
on   $l$-adic cohomology, which are  $\QQ_l$-linear bijections. The resulting polynomial
$$
P_1(C_0,T) = \det(1-\Phi T\mid H^1(C,\QQ_l))\in\QQ_l[T]
$$
is called the \emph{characteristic polynomial} of the curve $C$.
Recall that its degree equals the Betti number $b_1=\dim H^1(C,\QQ_l)$. Factoring $P_1(C,T)=\prod_{j=1}^{b_1}(1-\alpha_jT)$, one finds that the
reciprocal roots $\alpha_j\in\bar{\QQ}_l$ have absolute value $||\iota(\alpha_j)||=q^{m/2}$ 
for all complex embeddings $\iota:\QQ_l\ra\CC$.   
Note that extending the ground field by $\FF_q\subset\FF_{q^f}$ replaces $\Phi$ by the power $\Phi^f$ and
the reciprocal roots $\alpha_j$ by $\alpha_j^f$.
         
In general, it is  difficult to compute the characteristic polynomial. The following connects
characteristic polynomials with symmetries of the curve:

\begin{theorem}
\mylabel{characteristic polynomial}
Suppose there is a prime $l\neq p$ so that $H^1(C,\QQ_l)$ is irreducible as representation of $\Aut(C_0)$.
Then the characteristic polynomial  is of the form
$$
P_1(C_0,T)=(1\pm q^{1/2}T)^{b_1}.
$$
\end{theorem}

\proof
Set $A=\Aut(C_0)$. Choose a very ample invertible sheaf $\shL_0$ on $C_0$ endowed with
an $A$-linearizion, and   consider
the resulting embeddings $C_0\subset\PP^n_{\FF_q}$ and $A\subset\PGL_n(\FF_q)$.
The latter shows that $\Fr_{C_0}^m$ commutes with each automorphism $a\in A$.
Consequently, the induced bijection $\Phi:H^1(C,\QQ_l)\ra H^1(C,\QQ_l)$ is a morphism
of $A$-representation. By Schur's Lemma, it is multiplication by some scalar $\alpha^{-1}\in\QQ_l^\times$,
so the characteristic polynomial is $P_1(C,T)=(1-\alpha T)^{b_1}$, and $\alpha$ is a Weil $q$-number.
On the other hand, on knows that the characteristic polynomial has integral coefficients. Its linear coefficient
it $-b_1\alpha$, so $\alpha\in\QQ$. The only rational numbers that are also  Weil $q$-numbers    are $\pm q^{1/2}$.
\qed

\medskip
Recall that the curve $C$ is called \emph{supersingular} if its jacobian $J=\Pic^0_{C/k}$ is supersingular, that is, isogeneous
to a product of supersingular elliptic curve. We now have a criterion for supersingularity:

\begin{corollary}
\mylabel{irreducible supersingular}
Suppose there is a prime $l\neq p$ so that $H^1(C,\QQ_l)$ is irreducible as representation of $\Aut(C)$.
Then the curve $C$ is supersingular.
\end{corollary}

\proof
Replacing $\FF_q$ by a suitable finite extension, we may assume that $\Aut(C_0)=\Aut(C)$.
Then $\Phi=(1\pm q^{1/2}T)^{b_1}$ by the Theorem. As explained in \cite{Yui 1986},
the formal completion of the jacobian $J=\Pic^0_{C/k}$ is isogeneous to a selfproduct of the formal group scheme
$G_{1,1}$. The latter is equivalent to the supersingularity of $J$.
\qed

\medskip
Using Corollary \ref{cohomology irreducible}, we obtain as a special case:

\begin{corollary}
\mylabel{hermitian irreducible}
The Hermitian curves $C:\;y^q-y=x^{q+1}$ are supersingular.
\end{corollary}

%===========================================================
\section{Diagonal quotients}
\mylabel{Diagonal quotients}

We now collect some facts about quotients of products  of curves.
The situation is as follows: Fix an algebraically closed ground field $k$ of characteristic $p>0$,
and  let $C,C'$ be two smooth proper 
curves, say of genus $g,g'$, respectively. 
We   consider the   surface  $X=C\times C'$. The Betti numbers
of this smooth proper surface are
$$
b_1(X)=g+g'\quadand b_2(X)=2+gg'.
$$
Suppose that we have a finite group $G$, together with homomorphisms
$$
\Aut(C)\longleftarrow G\lra\Aut(C'),
$$ 
and consider the proper normal surface $Y=(C\times C')/G$. We assume   that
$G$ acts faithfully on each curve, hence freely on some open dense subsets.

\begin{proposition}
\mylabel{picard scheme}
Suppose that   the $G$-representations $H^1(C,\QQ_l)$ and $H^1(C',\QQ_l)$ 
contain  no copy of the trivial $G$-representation $\QQ_l$.
Then the Picard scheme $\Pic_{Y/k}$ of the proper normal surface $Y=G\backslash(C\times C')$ is zero-dimensional.
\end{proposition}

\proof
Since $Y$ is normal, the connected component $\Pic^0_{Y/k}$ is an extension of
an abelian variety by a zero-dimensional group scheme.
We thus have to show that $\Pic^0_{Y/k}$ contains no abelian variety.
Suppose to the contrary that it does. According to \cite{SGA 6}, Expos\'e XIII, Theorem 3.8, the pullback morphism
of group schemes $\Pic_{Y/k}\ra\Pic_{X/k}$ is quasiaffine. It follows that its image contains
an abelian variety. In turn, the homomorphism
$$
H^1(Y,\QQ_l)\lra H^1(X,\QQ_l)^G\subset H^1(X,\QQ_l)=H^1(C,\QQ_l)\oplus H^1(C',\QQ_l)
$$
is nonzero. By assumption, however, the $G$-invariant part of $H^1(X,\QQ_l)$ is trivial, contradiction.
\qed

\medskip
Recall that the \emph{N\'eron--Severi group} $\NS(Y)$ of a proper scheme $Y$ is the Picard group $\Pic(Y)$
modulo algebraic equivalence. This is a finitely generated abelian group, whose
rank is called the \emph{Picard number} $\rho(Y)$.

\begin{proposition}
\mylabel{picard number}
If $C,C'$ are supersingular, then   the 
Picard number   of the   proper normal surface $Y=G\backslash(C\times C')$ is given by $\rho(Y)=2+d$,
where  $d$ is the multiplicity of the trivial $G$-representation $\QQ_l$
in the tensor product $H^1(C,\QQ_l)\otimes H^1(C',\QQ_l)$.
\end{proposition}

\proof
Consider the following commutative diagram 
$$
\begin{CD}
\NS(X)\otimes\QQ_l @>>> H^2(X,\QQ_l)\\
@AAA @AAA\\
\NS(Y)\otimes\QQ_l @>>> H^2(Y,\QQ_l),
\end{CD}
$$
where the horizontal maps are induced from the $l$-adic cycle class maps. These maps are injective, 
as can be seen by using intersection numbers
on the smooth surface $X$ and the $\QQ$-factorial surface $Y$. Using the projection formula
for intersection numbers, one also sees that the vertical map on the left is injective.
According to  \cite{Yui 1986}, Theorem 2.5, the  upper horizontal map is bijective, because our curves $C,C'$ are supersingular.

It only remains to check that each each $G$-invariant class $[\shL]\in\NS(X)$ descends
to a class on $Y$, at least after passing to some multiple. Consider the commutative diagram with exact rows
$$
\begin{CD}
0 @>>> \Pic^0(X)^G @>>> \Pic(X)^G @>>> \NS(X)^G @>>> H^1(G,\Pic^0(X))\\
@.     @AAA                @AAA           @AAA \\
0 @>>> \Pic^0(Y)   @>>> \Pic(Y)   @>>> \NS(Y).
\end{CD}
$$
The group $\Pic^0(X)$ of algebraically trivial invertible sheaves is 
a divisible group, whence its cohomology vanishes. 
By a diagram chase, it suffices to check that each $G$-invariant invertible sheaf $\shL$ on $X$
descends to $Y$, at least after passing to some multiple.

To this end, consider the $G$-invariant open   subset $U\subset X$ on which the $G$-action is free,
and let $V\subset Y$ be its image. Clearly, the complements of these open subsets are zero-dimensional.
Since the action is free, we have a short exact sequence
$$
H^1(G,k^\times)\lra \Pic(V)\lra \Pic^G(U)\lra H^2(G,k^\times),
$$
and the outer terms are torsion groups. Passing to some multiple, we see that $\shL|_U$ descends to some invertible
sheaf $\shN_V$ on $V$.
Having only quotient singularities, the scheme $Y$ is $\QQ$-factorial. Passing again to some suitable multiple,
we see that $\shN_V$ extends to $Y$. This extension $\shN$ has the property that its preimage $\shN_X$ coincides
with $\shL$ over the open subset $U\subset X$ whose complement is zero-dimensional. It then follows that $\shN_X\simeq \shL$.
\qed

\medskip
The   projections $\pr_1:C\times C'\ra C$ induces a fibration 
$$
f:Y=(C\times C')/G\lra C/G.
$$
The generic fibers $Y_\eta=f^{-1}(\eta)$ is a smooth proper curve  over the function field $F=k(C/G)=k(C)^G$. 
By construction, $Y_\eta$ is a twisted form  of $C'\otimes_k F$. In particular, it is smooth of genus $g(Y_\eta)=g'$.
Combining the preceding two results, we obtain:

\begin{proposition}
\mylabel{generic fiber}
Suppose that the curves $C,C'$ are supersingular, and that the $G$-representations $H^1(C,\QQ_l)$ and $H^1(C',\QQ_l)$ 
contain  no copy of the trivial $G$-represen\-tation $\QQ_l$. Then the abelian group $\Pic^0(Y_\eta)$ is finitely generated of rank $d$,
where $d$ is the multiplicity of the trivial $G$-representation $\QQ_l$
in the tensor product $H^1(C,\QQ_l)\otimes H^1(C',\QQ_l)$.
\end{proposition}

\proof
The group $\Pic(Y)$ is finitely generated by Proposition \ref{picard scheme}. Its rank is   given by  $\rho=2+d$ by
Proposition \ref{picard number}.
Consider the restriction map $r:\Pic(Y)\ra\Pic(Y_\eta)$.  Its kernel is a finitely generated group of rank 
one, since all fibers of $Y\ra C/G$ are irreducible. Hence the image $\im(r)\subset\Pic(Y_\eta)$ is 
finitely generated of rank $1+d$. The cokernel of this map is $n$-torsion for some
integer $n\geq 1$, because the scheme $Y$ is $\QQ$-factorial. 
It remains to check that $\coker(r)$ is finitely generated. 
To this end, set $E=k(C) $, and regard $C_E$ as the
generic fiber of the projection $X=C\times C'\ra C$. This gives a commutative diagram
$$
\begin{CD}
@.       \Pic(Y) @>>> \Pic(Y_\eta) @>>> \coker(r) @>>> 0\\
@.       @VVV             @VVV              @VVV\\
0  @>>>   H  @>>> \Pic(X_E) @>>> \Pic(X_E)/H,
\end{CD}
$$
where $H$ is the image of the composite map $\Pic(Y)\ra\Pic(X)\ra\Pic(X_E)$.
The kernel of the vertical map in the middle is finite according to \cite{SGA 6}, Expos\'e XIII, Theorem 3.8.
By the Five Lemma, it thus suffices to check that the $n$-torsion in $\Pic(X_E)/H$
is finite. This follows from a general fact stated below.
\qed

\begin{lemma}
\mylabel{finitely generated}
Let $Z$ be a  proper, geometrically normal scheme over a ground field $E$, and $H\subset\Pic(Z)$ be a finitely generated subgroup.
Then the $n$-torsion subgroup of $\Pic(Z)/H$ is finite for all integers $n\geq 1$.
\end{lemma}

\proof
For any field extension $E\subset E'$, the induced morphism $\Pic(Z)\ra\Pic(Z\otimes E')$ is injective.
Thus we may assume that $E$ is algebraically closed. 
Set $D=\Pic(Z)$. The short exact sequence
$0\ra H\ra D\ra D/H\ra 0$ yields an exact sequence
$$
 \Tor_1(D,\ZZ/n\ZZ)\lra\Tor_1(D/H,\ZZ/n\ZZ)\lra H\otimes\ZZ/n\ZZ.
$$
Therefore, is suffices to check that the $n$-torsion in $D=\Pic(Z)$ is finite. This follows from
\cite{SGA 6}, Expos\'e XIII, Theorem 5.1, toether with  
the fact that the reduction of $\Pic^0_Z$ is an abelian variety.
\qed

%===========================================================
\section{Chern invariants}
\mylabel{Chern invariants}

In this final chapter we return to the relative minimal smooth model $\psi:S\ra\PP^1$
of the fibration $G\backslash(C\times C)\ra G\backslash C$, where $C:\, y^q-q=x^{q+1}$
is our Hermitian curve.
We saw in Corollary \ref{rational surface} that $S$ is a rational surface. Here
we compute its numerical invariants:

\begin{theorem}
\mylabel{chern invariants}
The Chern invariants for the rational surface $S$ are 
$$
c_2=e= q^2+q+6\quadand c_1^2=K_S^2=-q^2-q+6.
$$
The Picard number is given by the formula $\rho=q^2+q+4$.
\end{theorem}

\proof
We use Dolgachev's Formula \cite{Dolgachev 1972}
$$
e(S) = e(S_{\bar{\eta}})e(\PP^1) + e(S_\infty)-e(S_{\bar{\eta}})+\delta
$$
for the fibration $\psi:S\ra\PP^1$. Here $S_{\bar{\eta}}$ is the geometric generic fiber,
$S_\infty$ is the singular fiber, and $\delta=\delta(G,H^1(C,\FF_l)$ is the Swan conductor
for the $G$-representation $H^1(C,\FF_l)$ for some prime $l$ different from the characteristic $p$.
The curve $C$ has genus $q(q-1)/2$, thus  $e(S_{\bar{\eta}})=e(C)=-q^2+q+2$.
According to Proposition \ref{singular fiber}, the singular fiber is a tree of $q^2+4$ projective lines,
thus $e(S_\sigma)=q^2+5$. The Swan conductor is given by $\delta=q^2-1$ by \ref{wild conductor}.
Substituting these values into the right hand side of Dolgachev's Formula, we obtain
$e=q^2+q+6$. Since the surface $S$ is rational, we have $b_1=b_3=0$, thus $\rho=b_2=q^2+q+4$.

We may compute this   with the Tate--Shioda Formula for the fibration $S\ra\PP^1$ as well,
which gives $\rho=2+ (c-1) + r$, where $c$ is the number of irreducible components in the
singular fiber $S_\infty$, and $r$ is   the Mordell--Weil rank, that is, the rank
of $\Pic^0(S_\eta)$. We already saw that $c=q^2+4$, and  $r=q-1$ by Proposition \ref{generic fiber},
together with Corollary \ref{invariant cohomology}. Thus this second computation
also gives $\rho=q^2+q+4$.
\qed

%===========================================================

\end{document}